\documentclass[twocolumn]{IEEEtran}
\IEEEoverridecommandlockouts
\usepackage{xcolor}
\usepackage{booktabs, multirow, makecell}
\usepackage{bm, comment}
\usepackage[lofdepth,lotdepth]{subfig}
\usepackage{amsmath,amsfonts,amssymb}
\usepackage{mathabx,fixmath,algorithm, algorithmic}
\usepackage{pstool}
\usepackage{psfrag}
\usepackage{hyperref}
\usepackage{cite}
\usepackage{fancybox}
\usepackage{tikz}
\usetikzlibrary{shapes,arrows}
\usepackage{pgfplots}
\usepgfplotslibrary{groupplots}
\usetikzlibrary{calc}
\usepgfplotslibrary{statistics}
\usepackage{pgfplotstable}
\pgfplotsset{compat=1.8}

\def\x{\mathbold{x}}
\def\A{\mathbold{A}}

\def\bpi{\mathbold{\pi}}

\def\a{\mathbf{a}}
\def\b{\mathbold{b}}

\newcommand{\z}{\mathbold{z}}
\newcommand{\w}{\mathbold{w}}
\def\y{\mathbold{y}}
\def\e{\mathbold{e}}

\newcommand{\reals}{\mathbb{R}}

\DeclareMathOperator*{\argmin}{\textrm{argmin}}

\usepackage{theorem}

\newtheorem{theorem}{Theorem}

\newtheorem{assumption}{Assumption}

\newcommand{\rev}[1]{{\color{black} #1}}

\newcommand{\revi}[1]{{\color{black} #1}}

\newcommand{\cref}[1]{\eqref{#1}}

\usepackage{lipsum}
\usepackage{amsfonts}
\usepackage{graphicx}
\usepackage{epstopdf}
\usepackage{algorithmic}
\ifpdf
  \DeclareGraphicsExtensions{.eps,.pdf,.png,.jpg}
\else
  \DeclareGraphicsExtensions{.eps}
\fi

\renewcommand{\baselinestretch}{0.965}

\begin{document}
\title{Time-Varying Convex Optimization: Time-Structured Algorithms and Applications} 
%\title{Structured Time-Varying Convex Optimization for Large-Scale Dynamic Analytics and Applications} 
\author{Andrea Simonetto, Emiliano Dall'Anese, Santiago Paternain, Geert Leus, and Georgios B. Giannakis  \vspace*{-.9cm}

\thanks{

A. Simonetto is with the Optimization and Control Group of IBM Research Ireland, Dublin, Ireland. Email: andrea.simonetto@ibm.com. E. Dall'Anese is with the College of Engineering and Applied Science, University of Colorado, Boulder, CO 80309, USA. Email: emiliano.dallanese@colorado.edu. S. Paternain is with the Department of Electrical and Systems Engineering, University of Pennsylvania. Email: spater@seas.upenn.edu. G. Leus is with the Faculty of Electrical, Mathematics and Computer Science, Delft University of Technology, Delft 2628CD, The Netherlands. E-mail: g.j.t.leus@tudelft.nl. G.B. Giannakis is with the Digital Technology Center, University of Minnesota, Minneapolis, MN 55455 USA. Email: georgios@umn.edu.  }
}

\maketitle

% REQUIRED
\begin{abstract}
Optimization underpins many of the challenges that science and technology face on a daily basis. Recent years have witnessed a major shift from traditional optimization paradigms grounded on batch algorithms for medium-scale problems to challenging dynamic, time-varying, and even huge-size settings. This is driven by technological transformations that converted infrastructural and social platforms into complex and dynamic networked systems with even pervasive sensing and computing capabilities.
%Take smart grids for example, where continent-size networks are supposed to be optimized at fast temporal scales to capture the high variability of renewables. 
%Other paradigms emerge in autonomous systems, e.g., robotics, and intelligent transportation systems. From the analytics perspective, we have seen the drastic increase in the amount of data that is available, and its fast dynamic nature, which calls for novel optimization tools, of which time-varying optimization is one of paramount importance. 
%
The present paper reviews a broad class of state-of-the-art algorithms for time-varying optimization, with an eye to both algorithmic development and performance analysis. It offers a comprehensive overview of available tools and methods, and unveils open challenges in application domains of broad interest. The real-world examples presented include smart power systems, robotics, machine learning, and data analytics, highlighting domain-specific issues and solutions. The ultimate goal is to exempify wide engineering relevance of analytical tools and pertinent theoretical foundations.
\end{abstract}

\vspace*{-.5cm}

% !TEX root = revised_final.tex

%{\color{red} All: add any technical detail you think is reasonable to have

%REM AE: some of the technical/mathematical content could be expanded upon to make those parts accessible to non-specialists.}

\section{Introduction}\label{sec:introduction}

%%%%%%%%%%%%%%%%%%%%%%%%%%%%%%%%%%%%%%%%%%%%%%%%%%%%%%%%%%%%%%%%%%%%%%%%%%%%%%%%
%%%%%%%%%%%%%%%%%%%%%%%%%%%%%%%%%%%%%%%%%%%%%%%%%%%%%%%%%%%%%%%%%%%%%%%%%%%%%%%%
%%%%%%%%%%%%%%%%% INTRODUCTION 
%%%%%%%%%%%%%%%%%%%%%%%%%%%%%%%%%%%%%%%%%%%%%%%%%%%%%%%%%%%%%%%%%%%%%%%%%%%%%%%%
%%%%%%%%%%%%%%%%%%%%%%%%%%%%%%%%%%%%%%%%%%%%%%%%%%%%%%%%%%%%%%%%%%%%%%%%%%%%%%%%

%[Limit to three columns] \Andrea \Emiliano [all to help]

%\vskip3mm 

Optimization is prevalent across many engineering and science domains. Tools and algorithms from convex optimization have been traditionally utilized to support a gamut of data-processing, monitoring, and control tasks across areas as diverse as communication systems, power and transportation networks, medical and aerospace engineering, video surveillance, and robotics  -- just to mention a few. Recently, some of these areas -- and in particular  infrastructures such as power, transportation and communication networks, as well as social and e-commerce platforms -- are undergoing a foundational transformation, driven by major technological advances across various sectors, the information explosion propelled by online social media, and pervasive sensing and computing capabilities. Effectively, these infrastructures and platforms are revamped into complex systems operating in highly dynamic environments and with high volumes of heterogeneous information. This calls for revisiting several facets of workhorse optimization tools and methods under a different lens: the ability to process data streams and provide decision-making capabilities at time scales that match the dynamics of the underlying physical, social, and engineered systems using solutions that are grounded on conventional optimization methods can no longer be taken for granted. Take power grids, as a representative example:  economic optimization at the network level was performed using batch solvers at the minute or hour level to optimally dispatch large-scale fossil-fuel generation based on predictable loads; on the other hand, novel optimization tools are now desirable to carry network optimization tasks with solvers capable of coping with volatile renewable generation while managing the operation of a massive number of distributed energy resources. These considerations have spurred research and engineering efforts that are centered around \emph{time-varying optimization} -- a  formalism for modeling and solving optimization tasks in engineering and science under dynamic environments. 

%Optimization underpins many of the engineering challenges and problems that we face everyday. Recent years have witnessed a major shift from well-defined static and medium-scale optimization problems to online, time-varying, continent-size problems. Take smart grids for example, where continent-size networks are supposed to be optimized at fast temporal scales to capture the high variability of renewables~\cite{?}. Other paradigms emerge in autonomous systems, e.g., robotics~\cite{?}, and intelligent transportation systems~\cite{Alonso-Mora2017}. From the analytics perspective, we have seen the drastic increase in the amount of data that is available, and its fast dynamic nature, which calls for novel optimization tools, of which \emph{time-varying optimization} is one of paramount importance.  

Continuously-varying optimization problems represent a natural extension of time-invariant programs when the cost function and constraints may change continuously over time~\cite{Moreau1977,Polyak1987,Guddat1990,Dontchev2013}. Recently, time-varying optimization formalisms and the accompanying online solvers have been proposed both in continuous-time~\cite{Rahili2015,Fazlyab2018} and in discrete-time settings~\cite{Zavala2010,Simonetto2017a}. Their main goal is to develop algorithms that can track trajectories of the optimizers of the continuously-varying optimization program (up to asymptotic error bounds). The resultant algorithmic frameworks have demonstrated reliable performance in terms of convergence rates, with error bounds that relate tracking capabilities with  computational complexity; these features make time-varying algorithms an appealing candidate to tackle dynamic optimization tasks at scale, across many engineering and science domains. 

%Continuously varying optimization programs have appeared as a natural extension of time-invariant ones when the cost function, the constraints, or both, depend on a time parameter and change continuously in time~\cite{Moreau1977,Polyak1987,Guddat1990,Dontchev2013}. Recently, online algorithms have been put forward both in a continuous-time~\cite{Fazlyab2018,2} as well as in discrete-time settings~\cite{Zavala2010,Simonetto2017a} to track the trajectory of the optimizers of the continuously-varying optimization program, while it evolves in time up to asymptotic error bounds. These algorithms have demonstrated reliable performance in terms of convergence rates, error bounds, and computational effort; which makes them perfect candidates to tackle dynamic problems at scale. 

%%%%%%%%%%%%%%%%%%%%%%%%%%%%%%%%%%%%%%%%%%%%%%%%%%%%%%%%%%%%%%%%%%%%%%%%%%%%%%%%
%%%%%%%%%%%%%%%%%%%%%%%%%%%%%%%%%%%%%%%%%%%%%%%%%%%%%%%%%%%%%%%%%%%%%%%%%%%%%%%%
%%%%%%%%%%%%%%%%% FIGURE 1 
%%%%%%%%%%%%%%%%%%%%%%%%%%%%%%%%%%%%%%%%%%%%%%%%%%%%%%%%%%%%%%%%%%%%%%%%%%%%%%%%
%%%%%%%%%%%%%%%%%%%%%%%%%%%%%%%%%%%%%%%%%%%%%%%%%%%%%%%%%%%%%%%%%%%%%%%%%%%%%%%%

\begin{figure}
\centering
\includegraphics[width= 0.425\textwidth]{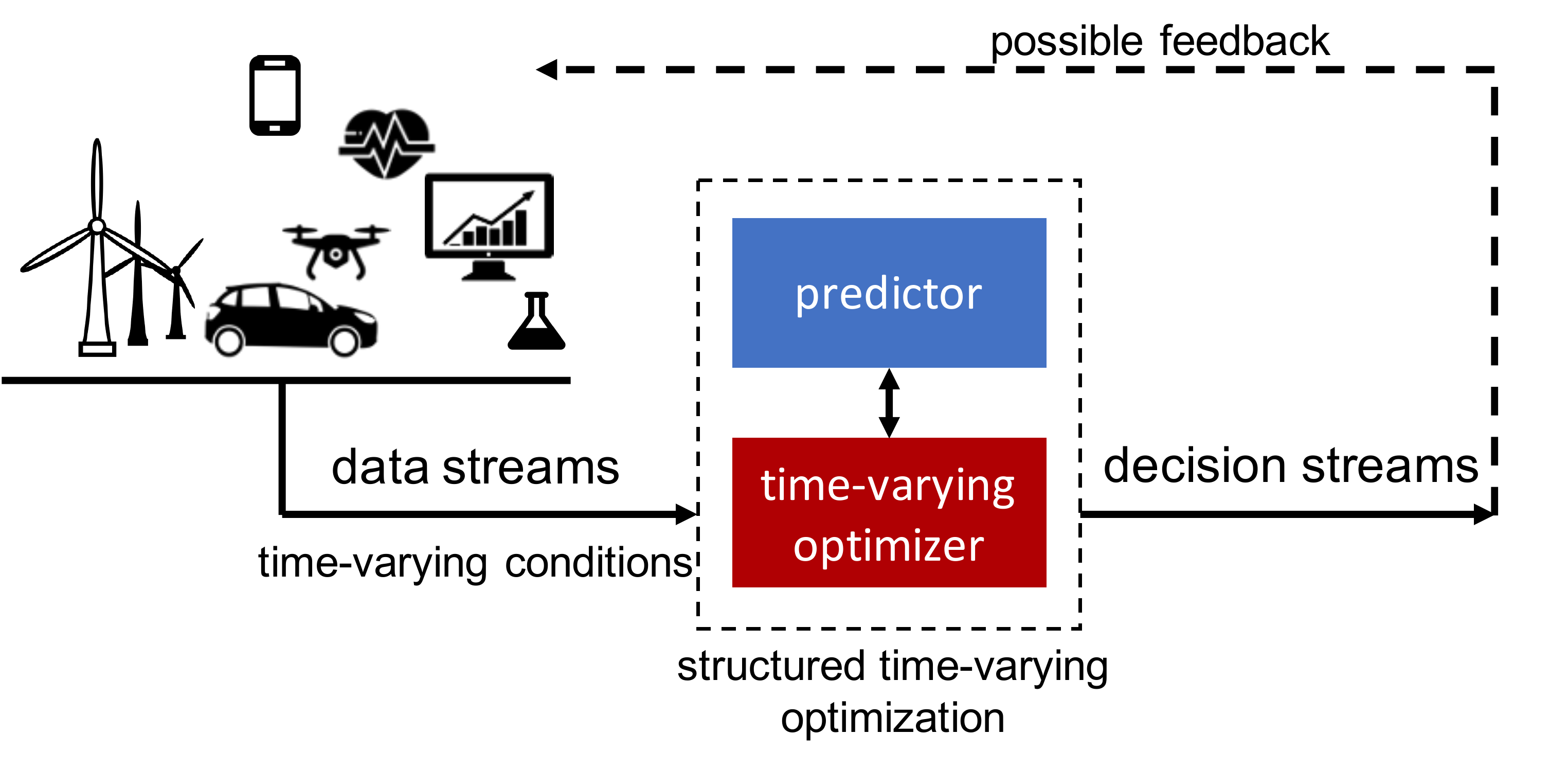}
\caption{\small The setup of time-varying optimization algorithms. Streaming data generated by time-varying systems are input to an optimizer. The optimizer can employ a predictor (that could be an oracle or a well-defined model), which feeds the optimizer with predictions of how the optimization problem will change. The optimizer then delivers a decision stream \rev{that is used to take actions that could be possibly fed back to affect the dynamical system operation.}}
\label{fig:depict}
\end{figure}

%%%%%%%%%%%%%%%%%%%%%%%%%%%%%%%%%%%%%%%%%%%%%%%%%%%%%%%%%%%%%%%%%%%%%%%%%%%%%%%%%%%%%%%%%%%%%%%%%%%%%%%%%%%%%%%%%%%%%%%%%%%%%%%%%%%%%%%%%%%%%%%%%%%%%%%%%%%%%%%%%%%%%%%%%%%%%%%%%%%%%%%%%%%

This paper overviews key modeling and algorithmic design concepts, with emphasis on \emph{time-structured} (structured for short) time-varying algorithms for convex time-varying optimization. The term ``structured'' here refers to algorithms that take advantage of the inherent temporal structure, meaning they leverage prior information (such as Lipschitz continuity or smoothness) on the evolution of the optimal trajectory  to enhance convergence and tracking. In contrast, the term ``unstructured'' will refer to time-varying algorithms that simply rely on current information  of cost and constraints. % and do not have means to leverage possible information regarding the temporal structure of the optimization problem. 
This also differentiates the present `time-structured' class from \emph{interactive} algorithms (that belong to the unstructured class), which are tailored to learner-environment or player-environment settings; e.g., the popular online convex optimization (OCO) setup~\cite{Shalev-Shwartz2012}, where online algorithms decide on current iterates (using only information of past cost functions), and subsequently the environment reveals partial  or  full information about the function to be optimized next. 

%%%%%%%%%%%%%%%%%%%%%%%%%%%%%%%%%%%%%%%%%%%%%%%%%%%%%%%%%%%%%%%%%%%%%%%%%%%%%%%%
%%%%%%%%%%%%%%%%%%%%%%%%%%%%%%%%%%%%%%%%%%%%%%%%%%%%%%%%%%%%%%%%%%%%%%%%%%%%%%%%
%%%%%%%%%%%%%%%%% FIGURE 2 
%%%%%%%%%%%%%%%%%%%%%%%%%%%%%%%%%%%%%%%%%%%%%%%%%%%%%%%%%%%%%%%%%%%%%%%%%%%%%%%%
%%%%%%%%%%%%%%%%%%%%%%%%%%%%%%%%%%%%%%%%%%%%%%%%%%%%%%%%%%%%%%%%%%%%%%%%%%%%%%%%
\begin{figure}
\centering
\includegraphics[width= 0.475\textwidth]{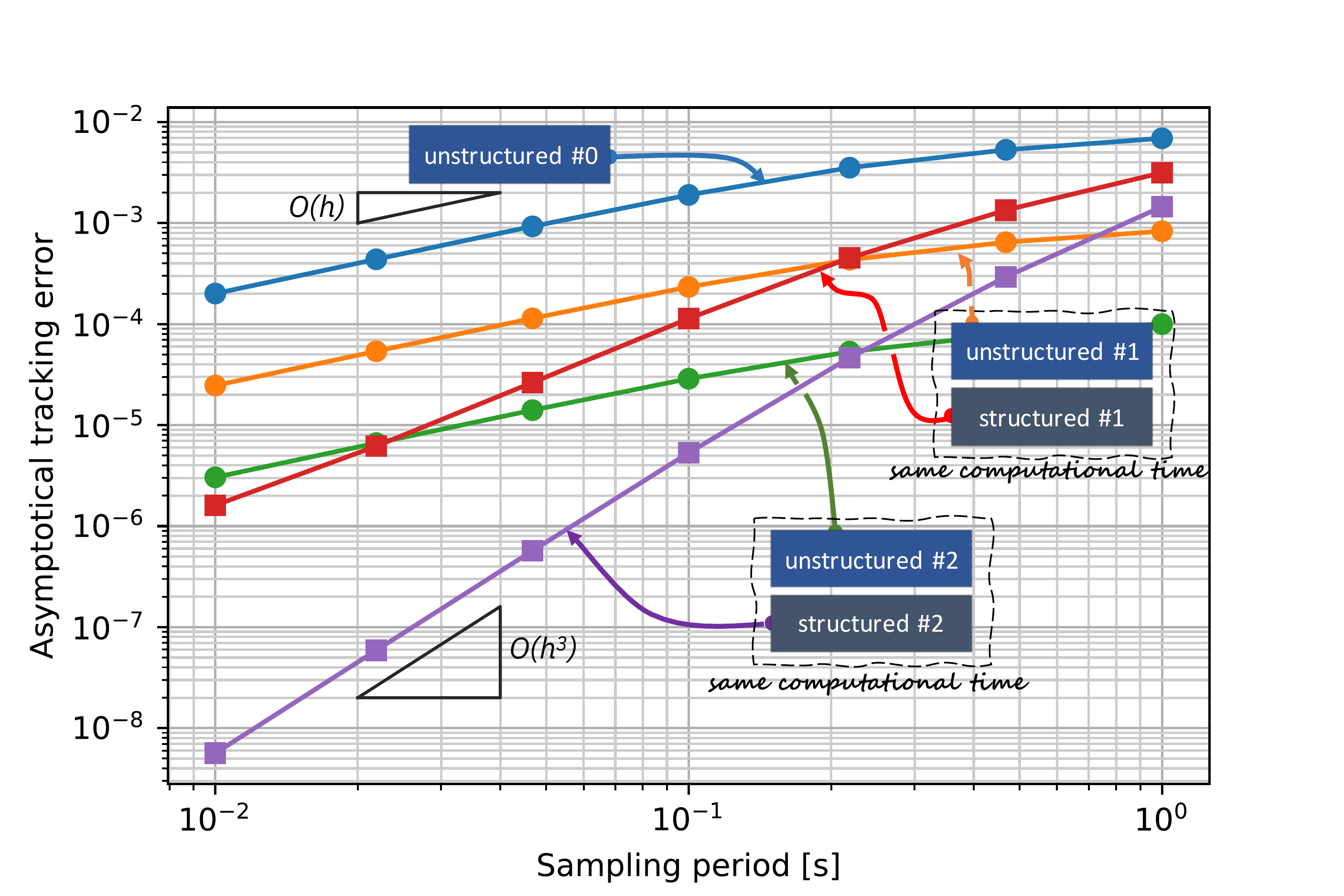}
\caption{\small Structured algorithms can outperform unstructured ones, even keeping the computational time fixed: here for a robot tracking problem. See text and footnote for description of the algorithms and \cite{Dixit2019,Bastianello2019} for the problem setting. Fig.~\ref{fig:example} will be referred to multiple times in the paper and the different elements will be clarified.}
\label{fig:example}
\end{figure}

%This paper overviews this online algorithms, with a particular emphasis on \emph{structured} time-varying algorithms for convex time-varying optimization. With structured we mean algorithms that make full use of the temporal structure of the problem, i.e., they exploit and/or model how the optimization problem evolves in time to better predict (anticipate) the optimizer trajectory. This in contrast with (more typical) unstructured time-varying algorithms, which only reacts to changes in a temporal-agnostic fashion. 

Figure~\ref{fig:depict} depicts a typical time-varying optimization setting. Streaming data are generated from time-varying systems, as in renewable generation that is intermittent, traffic conditions that change in transportation systems, or drop-off points for drone delivery that are mobile. Such settings inherit time variability in the optimization problem at hand. The optimizer can leverage a predictor (an oracle or a well-defined model), which feeds the optimizer with predictions of how the optimization problem may evolve over time. The optimizer then delivers a decision stream (i.e., an approximate optimizer) \rev{that is used to take operational actions such as committing a generator, or, adopting an optimal routing schedule for ridesharing vehicles. These actions could also affect and are therefore fed back} to the system \rev{(e.g., the optimal ridesharing schedule alters traffic and availability of vehicles in the future)}. When the input data streams are of large-scale and/or the decisions need to be made at a high frequency, traditional batch algorithms (that \revi{\emph{exactly}} solve the optimization problem at each time) are not viable because of underlying computational complexity bottlenecks. Hence, an online \emph{computationally frugal} optimization becomes essential to produce solutions in a timely fashion. 

\begin{figure}
\centering
\includegraphics[width= 0.475\textwidth]{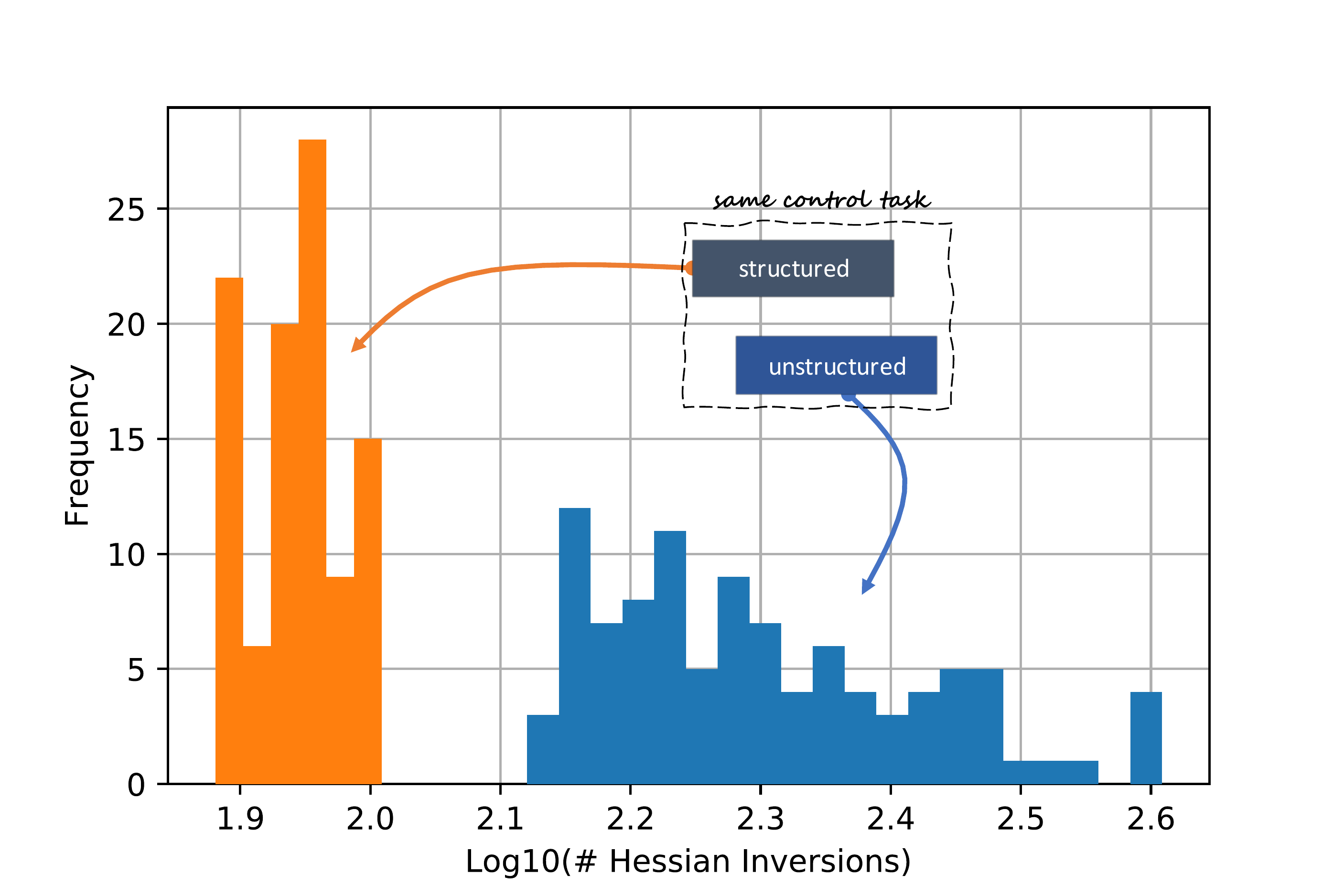}
  \caption{ \small Histograms with the total number of Hessian inversions required to control the Hicks reactor \cite{hicks1971approximation} for structured and unstructured MPC solvers. %The sampling time is 30s, the horizon is 10 and we run both controllers for 40 time steps. 
  Exploiting the temporal structure reduces the computational complexity, measured by the number of inversions of the Hessian, in solving the MPC. 
   }
                \label{fig_boxplot}
  \end{figure}

To further motivate structured time-varying methods, Figure~\ref{fig:example}\footnote{Unstructured algorithms 0, 1, and 2 are in this case online versions of the proximal gradient method, for which we perform 5, 7, and 9 passes of the methods, respectively. Structured algorithms here employ either a first- or a second-order Taylor model (for structured 1 and 2, respectively), and 5, and 20 passes of an online version of the proximal gradient method on a simplified quadratic problem;  see~\cite{Bastianello2019} for further details. } illustrates the asymptotic tracking error (asymptotic difference between optimal decisions and decisions delivered by some algorithms that will be described shortly) for different sampling periods ($h$) of discrete-time algorithms, for a robot tracking problem (see~\cite{Dixit2019} for the setting). The value of exploiting the temporal structure of the problem can be appreciated. Even keeping computational time fixed, structured algorithms outperform unstructured ones (here by several orders of magnitude). {\color{black}{Exploiting this structure may lead to a reduction of the computational cost of the algorithms. This is the case for instance when using model predictive control (MPC) on the Hicks reactor \cite{hicks1971approximation} (cf. Fig.~\ref{fig_boxplot} adapted from \cite{paternain2019prediction}).}}

%%%%%%%%%%%%%%%%%%%%%%%%%%%%%%%%%%%%%%%%%%%%%%%%%%%%%%%%%%%%%%%%%%%%%%%%%%%%%%%%
%%%%%%%%%%%%%%%%%%%%%%%%%%%%%%%%%%%%%%%%%%%%%%%%%%%%%%%%%%%%%%%%%%%%%%%%%%%%%%%% FIGURE 3
%%%%%%%%%%%%%%%%%%%%%%%%%%%%%%%%%%%%%%%%%%%%%%%%%%%%%%%%%%%%%%%%%%%%%%%%%%%%%%%%
%%%%%%%%%%%%%%%%%%%%%%%%%%%%%%%%%%%%%%%%%%%%%%%%%%%%%%%%%%%%%%%%%%%%%%%%%%%%%%%%

%%%%%%%%%%%%%%%%%%%%%%%%%%%%%%%%%%%%%%%%%%%%%%%%%%%%%%%%%%%%%%%%%%%%%%%%%%%%%%%%
%%%%%%%%%%%%%%%%%%%%%%%%%%%%%%%%%%%%%%%%%%%%%%%%%%%%%%%%%%%%%%%%%%%%%%%%%%%%%%%% 
%%%%%%%%%%%%%%%%%%%%%%%%%%%%%%%%%%%%%%%%%%%%%%%%%%%%%%%%%%%%%%%%%%%%%%%%%%%%%%%%
%%%%%%%%%%%%%%%%%%%%%%%%%%%%%%%%%%%%%%%%%%%%%%%%%%%%%%%%%%%%%%%%%%%%%%%%%%%%%%%%

The main goal of this overview paper is threefold: 
\begin{itemize}
\item[\emph(i)] To expose models and algorithms for structured time-varying optimization settings, from both analytical and an application-oriented perspectives; 

\item[\emph(ii)] Demonstrate applications of structured time-varying optimization algorithms (and deep dive into two, namely a robotic and a power system application); and

\item[\emph(iii)] Draw links with the growing landscape of unstructured algorithms for dynamic optimization problems.

\end{itemize}

%The overarching objective is to present how structured time-varying optimization is essential and could make a critical difference in many applications

\smallskip
\rev{
%%%
%%% NOTATION
\noindent {\bf \revi{Setting and notation.}} We deal with convex optimization problems~\cite{Boyd2004a,Rockafellar1970}, as well as first-order algorithms~\cite{Nesterov2004}. %Our problems live generally in $\reals^n$. 
Vectors are represented with $\x\in\reals^n$, and the Euclidean norm is indicated as $\|\cdot\|$. \revi{We mainly deal with strongly convex and smooth functions}. A function $f: \reals^n \to \reals$ is $m$-strongly convex for a constant $m>0$, i.e., $f(\x)-m/2\|\x\|^2$ is convex, and $L$-smooth for a constant $L>0$ iff its gradient is $L$-Lipschitz continuous or equivalently iff $f(\x)-L/2\|\x\|^2$ is concave. Sometimes, we deal with extended-real-valued functions $\varphi: \reals^n \to \reals\cup\{+\infty\}$ (which can explicitly admit infinite values, e.g., the indicator function). We define the \textit{subdifferential} of $\varphi$ as the set
$
	\x \mapsto \left\{ \z \in \reals^n\ |\ \forall \y \in \reals^n:\ \langle \y-\x, \z \rangle + \varphi(\x) \leq \varphi(\y) \right\}.
$ 
Given a convex set $\mathcal{X}$, $\mathrm{proj}_{\mathcal{X}}\{\x\}$ denotes a closest point to $\x$ in $\mathcal{X}$, namely $\mathrm{proj}_{\mathcal{X}}\{\x\} \in \arg \min_{\y \in \mathcal{X}} \|\x - \y\| $. We also use $O(\cdot)$ to represent the big-$O$ notation.

%%%

}

% !TEX root = revised_final.tex

\section{Time-varying optimization}\label{sec:tvo}

%%%%%%%%%%%%%%%%%%%%%%%%%%%%%%%%%%%%%%%%%%%%%%%%%%%%%%%%%%%%%%%%%%%%%%%%%%%%%%%%
%%%%%%%%%%%%%%%%%%%%%%%%%%%%%%%%%%%%%%%%%%%%%%%%%%%%%%%%%%%%%%%%%%%%%%%%%%%%%%%%
%%%%%%%%%%%%%%%%% TIME-VARYING OPTIMIZATION
%%%%%%%%%%%%%%%%%%%%%%%%%%%%%%%%%%%%%%%%%%%%%%%%%%%%%%%%%%%%%%%%%%%%%%%%%%%%%%%%
%%%%%%%%%%%%%%%%%%%%%%%%%%%%%%%%%%%%%%%%%%%%%%%%%%%%%%%%%%%%%%%%%%%%%%%%%%%%%%%%

%\textbf{[3 pages]} \Andrea \Emiliano \Santiago [all to help]

%\vskip3mm

Let $f: \reals^n \times \reals_{+} \to \reals$ be a convex function parametrized over time, i.e., $f(\x; t)$, where $\x \in \reals^n$ is the decision variable and $t\geq 0$ is time. Let $X(t)\subseteq \reals^n$ be a convex set, which may also change over time. We are interested here in solving:% the following:
\begin{equation}\label{tvp}
\min_{\x\in X(t)} \, f(\x; t), \quad \textrm{for all } t\geq 0.
\end{equation}
To simplify exposition, we assume that the cost function $f$ is $m$-strongly convex for all $t$ (this is nevertheless a standard assumption in most prior works), and that the constraint set is never empty. With these assumptions in place, at any time $t$, Problem~\eqref{tvp} has a unique global \rev{optimizer}. This translates to finding the optimal solution trajectory  
\begin{equation}\label{eq:problem}
\x^{\star}(t) := \argmin_{\x\in X(t)} \, f(\x; t), \quad \textrm{for all } t\geq 0.
\end{equation}

As an example, for the robot tracking problem for which the results have been shown in Figure~\ref{fig:example}, $f(\x;t)$ is a time-varying performance metric for the tracking performance of a robot formation that is following a robot leader; for example, $f(\x;t) = \|\x - \b(t)\|^2 + \mathcal{R}(\x)$, where $\mathcal{R}(\x)$ is some pertinent regularization function and $\b(t)$ encodes the tracking signal. On the other hand, $X(t)$ represents some physical or hardware constraints for the robots. At each $t'$, the information available is $\{f(\x;t), t \leq t' \}$ and $\{X(t), t \leq t'\}$; based on a possibly limited computational complexity, and without any information regarding future costs and constraints, the next decision $\x(t')$ has to be made; the objective is to produce a decision $\x(t')$ that is as close as possible to $\x^{\star}(t')$. 

If Problem~\eqref{eq:problem} changes slowly, and sufficient computational power is available,  existing batch optimization methods may identify the optimal trajectory $\x^{\star}(t)$; for example, if the parameter $\b(t)$ above exhibits step changes every 10 seconds, and a distributed batch algorithm converges in 5 seconds, then $\x^{\star}(t)$ can be identified (within a given accuracy).  On the other hand, in highly dynamic settings, computational and communication bottlenecks may prevent batch methods to produce solutions in a timely manner (e.g., $\b(t)$ changes every 0.5 seconds, and a distributed batch algorithm converges in 5 seconds); the problem then becomes related to the synthesis of computationally-affordable algorithms that can produce an approximate optimizer trajectory $\hat{\x}(t)$ on the fly; accordingly, a key performance of these algorithms is the ``distance''  between the approximate solution trajectory $\hat{\x}(t)$ and the optimal one $\x^{\star}(t)$. 

\rev{
\smallskip %%%%%%%%%%%%%%%%%%%%%%%%%%%%%%%%%%%%%%%%%%%%%%%%%%%

\noindent {\bf Time-structured and time-unstructured algorithms.} The term ``structured'' refers to algorithms that, at time $t$, exploit a (learned) model to \emph{predict} how the optimizer trajectory $\hat{\x}(t)$ evolves, say from $t$ to $t'$, and then \emph{correct} the prediction by approximately solving the optimization problem \emph{obtained} at $t'$. Unstructured algorithms instead have no evolution model and use only the optimization problems that are revealed at each time.  
A useful parallelism is the Kalman filter versus the recursive least-squares (RLS) estimator. While the Kalman filter is endowed with a model to predict how the state evolves in time, and then corrects the prediction with new up-to-date observations, RLS relies solely on the observations. Structured time-varying algorithms leverage an evolution model to predict and observe new problems to correct their predictions. Unstructured ones rely only on observations.    
}

%At every time $t'$, one uses the information that they have available up to that time (including the function $f(\x; t')$) to generate an approximate optimizer $\hat{\x}(t')$. When the optimization problems at fixed $t$ are of large-scale and/or we need approximate optimizers at high frequency, it is unlikely that one can solve the fixed-time problem at time $t'$ at optimality, and therefore it is essential to build online algorithms that can converge and track the optimizer trajectory $\x^{\star}(t)$ while it evolves in time. 

\smallskip %%%%%%%%%%%%%%%%%%%%%%%%%%%%%%%%%%%%%%%%%%%%%%%%%%%

\noindent {\bf Performance metrics.} Different performance metrics can be considered for  online algorithms that generate approximate trajectories for Problem~\eqref{eq:problem}. \rev{They all capture the fact that the computation of $\hat{\x}(t)$ is time-limited, computationally-limited, or both, and therefore $\hat{\x}(t)$ is an approximate optimizer at time $t$. Here, it is more fruitful to look at the computation of $\hat{\x}(t)$ as limited by time: to compute $\hat{\x}(t)$ one has at most $\Delta t$.}

An immediate \rev{performance metric} is the \emph{asymptotical tracking error} (ATE), defined as
\begin{equation}\label{eq_ate}
\textrm{ATE} : = \limsup_{t\to\infty} \|\hat{\x}(t) - \x^{\star}(t)\|,
\end{equation}
which captures how the algorithm performs in an asymptotic sense. In general, one seeks \emph{asymptotic consistency} of the algorithm, i.e., if $\x^{\star}(t)$ is asymptotically stationary, then the ATE should be zero. However, if $\x^{\star}(t)$ is time-varying, the ATE cannot be zero for unstructured algorithms, while it could be zero for structured algorithms\footnote{A dynamic regret notion based on the cost function is also available, but we do not discuss this here. The interested reader is referred to e.g.~\cite{Shahrampour2018,Bernstein2019}.}.

A second metric that is relevant for time-varying optimization problems is the \emph{time rate} (TR), defined as
\begin{equation}
\textrm{TR} : = \frac{\textrm{time required for the computation of $\hat{\x}(t)$}}{\textrm{time allowed for the computation of $\hat{\x}(t)$}}.
\end{equation}
\rev{Here we define as ``time required,'' the time needed for the computation of an approximate $\hat{\x}(t)$, which delivers a predefined ATE}. The TR is a key differentiator for time-varying optimization: online algorithms need to be able to deliver an approximate $\hat{\x}(t)$ in the allocated time. Data streams generate decision streams with the same frequency, and the online optimization algorithm needs to have a TR less than one to be implementable. \rev{The TR sets also an important trade-off between ATE and implementability. One typically cannot expect a very low ATE and implementable solutions. }

The third metric is the \emph{convergence rate} (CR), which can be informally defined as
\begin{equation}
\textrm{CR} : = \textrm{``how fast'' an algorithm converges to the ATE}.
\end{equation}
Convergence rate will be formalized for discrete-time algorithms and continuous-time algorithms shortly. For discrete-time algorithms, under current modeling  assumptions, it will be possible to derive  $Q$-linear convergence results (definition given later on); on the other hand, for  continuous-time algorithms the convergence rate will be exponential and related to the exponent of a carefully constructed Lyapunov function. 

Typically, the algorithmic design will involve a trade-off between the ATE and CR; for instance,  lower levels of ATE may be achievable at the expense of a higher CR. CR is then important, not only at the start, but also when abrupt changes happen (and then the CR captures how fast the algorithm responds to those changes and disturbances). 

An additional metric is a measure that distinguishes between structured and unstructured algorithms, here referred to as structure gain (SG). \revi{It could be defined as} the ratio between the ATE obtained with a structured method divided by the ATE obtained with a competing unstructured method; that is: 
\begin{equation}
\textrm{SG} : = \frac{\textrm{ATE for \revi{selected} structured method}}{\textrm{ATE for \revi{competing} unstructured method}}.
\end{equation}
Of course, both algorithms are  constrained to use the same computational time for  $\hat{\x}(t)$. This metric assists in the decision as to whether  to use \revi{the selected} structured or \revi{the competing} unstructured algorithm for a given time-varying optimization  task. 
\begin{figure*}
\centering
\includegraphics[width=0.95\textwidth]{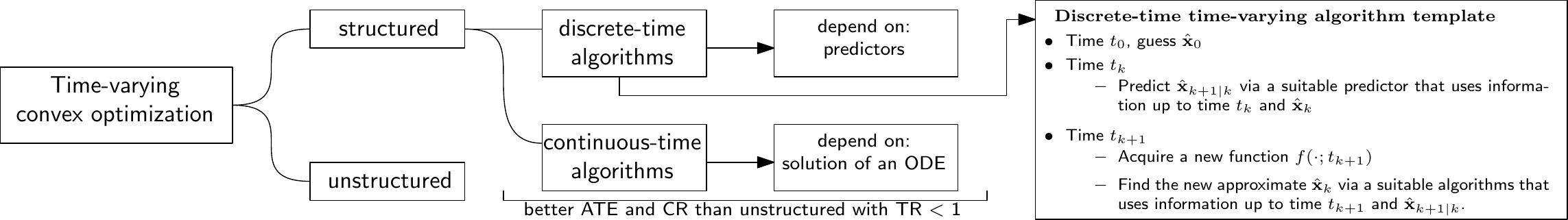}
\caption{\small{Algorithms presented in this paper.}}
\label{fig-3-str}
\end{figure*}
We have already seen in Figure~\ref{fig:example} that the value of structure can lead to an SG greater than one, further motivating the use of structured methods. 

In Figure~\ref{fig-3-str}, a general overview of the algorithms that will be presented in this paper is given together with their connections.

%%%%%%%%%%%%%%%%%%%%%%%%%%%%%%%%%%%%%%%%%%%%%%%%%%%
%%% SUBSECTION: DISCRETE-TIME ALGORITHMS
%%%%%%%%%%%%%%%%%%%%%%%%%%%%%%%%%%%%%%%%%%%%%%%%%%%

\subsection{Discrete-time algorithms} 

This section surveys \emph{discrete-time} algorithms. %Most of it is also valid for continuous-time algorithms, but somewhat more natural if presented in discrete-time. 
Consider sampling Problem~\eqref{eq:problem} at defined sampling times $\{t_k = k h, k \in \mathbb{N}\}$, with $h$ the sampling period; thus, one arrives at a sequence of time-invariant problems:
\begin{equation}\label{eq:tip}
\x^{\star}(t_k) := \argmin_{\x\in X} \, f(\x; t_k), \,\,\,\, t_k = k h, k \in \mathbb{N} \, .
\end{equation}
For simplicity of exposition, we drop the time dependency of the constraints and consider static sets. As long as one can solve each (time-invariant)  Problem~\eqref{eq:tip} within an interval $h$ using existing algorithms, then a ``batch solution mode'' is sufficient to identify the optimal trajectory $\{\x^{\star}(t_k), k \in \mathbb{N}\}$. This batch approach is, however, hardly viable, except for low-dimensional problems that can be sampled with sufficiently long sampling periods (i.e., when the problem changes sufficiently slowly). We focus here on the case where  one can afford only one or a few steps of a given algorithm within an interval $h$ -- i.e., an online approach. This setting can then be cast as the problem of synthesizing online algorithms that can \emph{track} $\{\x^{\star}(t_k), k \in \mathbb{N}\}$, within a given ATE.

%We will pursue an online approach, which will find approximate solutions of each of the time-invariant problems~\eqref{eq:tip} and eventually will get close to the optimizer trajectory. 

A key assumption for any online approach is that the difference between solutions at two consecutive times  is \emph{bounded}:
\begin{assumption}\label{as:1}
The distance between optimizers at subsequent times is uniformly upper bounded as:
$$
\|\x^{\star}(t_k) - \x^{\star}(t_{k-1})\| \leq K, \quad \forall k>0, \, K< \infty.
$$
\end{assumption}
The constant $K$ will play a key role in the ATE, as shown shortly. Assumption~\ref{as:1} is general, inasmuch it does not forbid the underlying trajectory $\x^{\star}(t)$ to have finite jumps\footnote{
\rev{
Meaning that $\x^{\star}(t)$ can be discontinuous in time, but the discontinuity has to be bounded, so that As.~\ref{as:1} holds for the choice of sampling period.  
}
}. 

A stronger assumption, often required in time-structured optimization, is that the time derivative of the gradient of the cost function\footnote{This can be generalized for a non-smooth cost function of the form $f(\x;t) + g(\x)$, as long as $f(\x;t)$ is differentiable, e.g., $\|x-t\|^2 + |x|$~\cite{Bastianello2019}.}, i.e., $\nabla_{t \x}f(\x;t)$, is bounded.
\begin{assumption}\label{as:2}
%The following bound is in place 
For all $t$ and all $\x$:
$
\|\nabla_{t \x}f(\x;t)\| \leq \Delta_0 < \infty.
$
\end{assumption}
Assumption~\ref{as:2}, along with $m$-strong convexity of the cost function, guarantees that the trajectory $\x^{\star}(t)$ is globally Lipschitz in time~\cite{Dontchev2009, Bastianello2020}, and in particular
\begin{equation}
\|\x^{\star}(t') - \x^{\star}(t)\| \leq \frac{\Delta_0}{m} |t'-t| . 
\end{equation}
This is key for structured time-varying algorithms, \rev{and typically not required in unstructured algorithms or in online convex optimization (OCO)~\cite{Shalev-Shwartz2012}}. Notice further that As.~\ref{as:2} implies As.~\ref{as:1} with the choice $K = \Delta_0\, h/m$. 

In this discrete-time setting, an online algorithm will generate a sequence of approximate optimizers. 
Hereafter, we will denote as $\hat{\x}_k$ the output of the algorithm at time $t_k$ for simplicity, while the sequence is denoted as $(\hat{\x}_k)_{k\in\mathbb{N}_{+}}$. Different algorithms will be distinguished based on which predictor they use and how they generate $\hat{\x}_{k}$. 

%\begin{figure}

%%%%%%%%%%%%%%%%%%%%%%%%%%%%%%%%%%%%%%%%%%%%%%%%%%%%%%%%%%%%%%%%%%%%%%%%%%%%%%
%%%%%%%%%%%%%%%%%%% ALGORITHM
%%%%%%%%%%%%%%%%%%%%%%%%%%%%%%%%%%%%%%%%%%%%%%%%%%%%%%%%%%%%%%%%%%%%%%%%%%%%%%%
%\scriptsize
%\centering
%\noindent\ovalbox{\begin{minipage}{0.43\textwidth}
%\centering
%\textsf{Discrete-time time-varying algorithm template\, \color{red} Make table or remove?} 
%\end{minipage}}
%\cornersize{.2}

%\noindent\ovalbox{\begin{minipage}{0.43\textwidth}
%\begin{itemize}
%\item Time $t_0$, guess $\hat{\x}_0$
%\item Time $t_{k}$
%\begin{itemize}
%\item Predict $\hat{\x}_{k+1|k}$ via a suitable predictor that uses information up to time $t_k$ and $\hat{\x}_k$
%\end{itemize}
%\item Time $t_{k+1}$
%\begin{itemize}
%\item Acquire a new function $f(\cdot; t_{k+1})$
%\item Find the new approximate $\hat{\x}_{k}$ via a suitable algorithms that uses information up to time $t_{k+1}$ and $\hat{\x}_{k+1|k}$. 
%\end{itemize}
%\end{itemize}
%\end{minipage}
%}
%\normalsize

%\vspace*{-.5cm}
%\end{figure}
%\vskip2mm

%%%%%%%%%%%%%%%%%%%%%%%%%%%%%%%%%%%%%%%%%%%%%%%%%%%%%%%%%%%%%%%%%%%%%%%%%%%%%%%

\smallskip

\noindent {\bf No-predictor algorithms.} In this case, online algorithms do not have a  ``prediction'' step; rather, they only perform ``corrective'' steps once the cost function is acquired. These algorithms are called in different ways (among which catching up, running, correction-only, unstructured) and probably firstly appeared %in the context of optimization 
with Moreau~\cite{Moreau1977}.  For example, a running projected gradient to approximately solve~\eqref{eq:tip} is given by the recursion 
\rev{
%\begin{subequations}
\begin{equation}\label{prj.gr}
\hat{\x}_{0} \!=\! {\bf 0}, \quad 
\hat{\x}_{k} \!=\! \textrm{proj}_{X}\{\hat{\x}_{k-1} \!-\! \alpha \nabla_{\x} f(\hat{\x}_{k\!-\!1}; t_k)\}, \, k \in \mathbb{N},
\end{equation}
%\end{subequations}
where} $\textrm{proj}_{X}\{\cdot\}$ denotes the projection operator and $\alpha$ is a carefully chosen step size \revi{(that could be time-varying as well)}. \revi{In~\eqref{prj.gr}, the projected gradient is applied one time per time step $k$, but one could also apply multiple gradient steps, say $C$, per time step. Notwithstanding this, } in general, these unstructured discrete-time algorithms achieve a high ATE. To formalize this result, we focus on a class of algorithms that exhibit a  Q-linear convergence. In particular, let $\mathcal{M}$ be an algorithm that when applied to $\hat{\x}_{k}$ at time $t_{k+1}$ for function $f(\x; t_{k+1})$ produces an $\hat{\x}_{k+1}$ for which, 
\begin{equation}\label{eq:linear}
\|\hat{\x}_{k+1} - \x^{\star}(t_{k+1})\| \leq \varrho \|\hat{\x}_{k} - \x^{\star}(t_{k+1})\|, \quad \varrho \in (0,1); 
\end{equation}
then algorithm $\mathcal{M}$ is called Q-linear convergent. This class is common in time-varying optimization (e.g., projected gradient\rev{~\eqref{prj.gr}} is Q-linear on a $m$-strongly convex,  $L$-smooth cost function~\cite{Nesterov2004}). \revi{When the algorithm $\mathcal{M}$ is then applied $C$ times (as e.g., in~\eqref{eq:sg-corr}), we obtain: $\|\hat{\x}_{k+1} - \x^{\star}(t_{k+1})\| \leq \varrho^C \|\hat{\x}_{k} - \x^{\star}(t_{k+1})\|$ }. 
%,Ryu2015,Taylor2017}). 
The following general result is in place. 

\begin{theorem}\label{th.1}\emph{(Informal)} Let $\mathcal{M}$ be an optimization algorithm that converges Q-linearly as in~\eqref{eq:linear}. Then, under Assumption~\ref{as:1}, the same algorithm $\mathcal{M}$ applied $C$ times for each time $t_{k}$, converges Q-linearly to the optimizer trajectory of a time-varying problem up to an error bound as
$$
\|\hat{\x}_{k+1} - \x^{\star}(t_{k+1})\| \leq \varrho^C( \|\hat{\x}_{k} - \x^{\star}(t_{k})\| + K), 
$$
and
$
\limsup_{k\to\infty} \|\hat{\x}_{k} - \x^{\star}(t_{k})\| = \varrho^C O(K) = \frac{\Delta_0}{m} \varrho^C O(h)
,
$
where the last equality is valid under Assumption~\ref{as:2}.
\end{theorem}

\rev{
\begin{proof}
\emph{(Sketch)} At time $t_k$, if algorithm $\mathcal{M}$ is applied $C$ times, starting on $\hat{\x}_{k}$ and ending at $\hat{\x}_{k+1}$, by Q-linear convergence of $\mathcal{M}$, we can write
\begin{multline*}
\|\hat{\x}_{k+1} - \x^{\star}(t_{k+1})\| \leq \varrho^C( \|\hat{\x}_{k} - \x^{\star}(t_{k+1})\|) \leq \\ \varrho^C( \|\hat{\x}_{k} - \x^{\star}(t_{k})\| + \|\x^{\star}(t_{k+1}) -\x^{\star}(t_{k})\|), 
\end{multline*}
and by using Assumption~\ref{as:1} the first claim is established. The second claim is proved by recursively applying the first claim, and by geometric series summation. 
\end{proof}
}

\smallskip

The results of the theorem are general and assert that the sequence $(\hat{\x}_{k})$ tracks the solution trajectory up to a ball of size $\varrho^C O(K)$. If $C \to \infty$, the time-invariant problem is solved exactly and we are back to the batch mode (and the error is $0$), i.e., the time-varying algorithm is asymptotically consistent. If As.~\ref{as:2} holds true, then the asymptotic error is proportional to the sampling period $h$ (cf. Figure~\ref{fig:example}). \rev{In addition, \emph{for fixed $\varrho \in (0,1)$, $C<\infty$, and if the path-length
$
\sum_{k=1}^T\!\|\x^{\star}(t_{k})\!-\!\x^{\star}(t_{k-1})\|
$
grows at least linearly in $T$, no unstructured method of this type can reach a zero  
ATE}~\cite{Besbes2013, Li2018a}.
}

\smallskip

\noindent{\bf Predictors. } We now focus on discrete-time algorithms that are endowed with a prediction. Various predictors are considered, \revi{and we will call as $\hat{\x}_{k+1|k}$ the predicted decision variable for time $t_{k+1}$ with information up to time $t_k$.}  

\smallskip
\noindent\emph{Clairvoyant oracles and expert oracles.} Clairvoyant oracles offer an exact prediction: i.e., they provide a $\hat{\x}_{k+1|k}$, for which $ \|\hat{\x}_{k} - \x^{\star}(t_{k}) \| = \|\hat{\x}_{k+1|k} - \x^{\star}(t_{k+1})\|$, \rev{as if they knew the function $f(\cdot; t_{k+1})$ and its gradient at time $t_k$}. In this context, clairvoyant oracles completely remove the time effect in the optimizer and the optimizer can proceed as if the cost function were not varying in time. Clairvoyant oracles are impractical (they need to have a perfect knowledge of the future), \revi{but they offer good performance lower bounds (since, one cannot do better than them)}. A noteworthy example of when one can use a clairvoyant oracle is when the cost function has a time drift, i.e., $f(\x; t) = f(\x+ \mathbold{\alpha}t)$, and the oracle can estimate the drift vector $\mathbold{\alpha}$ exactly based on historical data. 

Expert oracles, \rev{hints, or predictable sequences are considered, e.g., in~\cite{Rakhlin2013, Jadbabaie2015, Dekel2017}. In~\cite{Jadbabaie2015},} one has access to a sequence $(\mathbold{m}_k)_{k \in \mathbb{N}_+}$ of gradient approximators. When $\mathbold{m}_k = {\bf 0}$, \rev{i.e., meaning no knowledge or prediction about the future}, we recover an unstructured algorithm. When $\mathbold{m}_{k} = \nabla_{\x}f(\x; t_{k})$ at time $t_k$, then one recovers the online algorithm of~\cite{Chiang2012}. Finally, when $\mathbold{m}_{k} = \nabla_{\x}f(\x; t_{k+1})$, one recovers a clairvoyant oracle. Based on the error $\|\mathbold{m}_{k} - \nabla_{\x}f(\x; t_{k+1})\|$, one can then derive dynamic ATE results.

\smallskip

\noindent\emph{Model-based predictors.} These predictors are built on a model of the variations of the cost function, or of its parameters. 

\smallskip

$\bullet$ Prediction based on first-order optimality conditions~\cite{Dontchev2013, Zavala2010, Paper1, Simonetto2017a, Paper4, Qi2019}. A large class of predictors comes from deriving models based on first-order optimality conditions. We could call these predictors environment-agnostic, since they are not interested in modeling how the environment changes, but only how the optimization problem is affected. To introduce these predictors, let us consider an unconstrained problem (easier than Problem~\eqref{eq:tip}) as:
\begin{equation}\label{eq:un}
\x^{\star}(t_{k}) = \argmin_{\x\in\reals^n} f(\x;t_{k}). 
\end{equation}
To derive a model for how the problem is changing from $t_k$ to $t_{k+1}$, we look at the first-order optimality conditions at time $t_{k}$, which can be framed as
\begin{equation}\label{eq:fermat}
	\nabla_{\x} f(\x;t_{k}) = {\bf 0}.
\end{equation}
To predict, how this first-order optimality condition changes in time, with information available up to $t_k$, we use a Taylor expansion around $(\hat{\x}_k; t_k)$ as 
\begin{align}\label{eq:split}
\begin{split}
	{\bf 0}  = & \nabla_{\x} f(\x;t_{k+1}) \approx \varphi_k(\x) := \nabla_{\x} f(\hat{\x}_k;t_k) + \\ & + \nabla_{\x\x} f(\hat{\x}_k;t_k) (\x - \hat{\x}_k) + h\,  \nabla_{t \x} f(\hat{\x}_k;t_k),
\end{split}
\end{align}
where it is assumed that the Hessian $\nabla_{\x\x} f(\hat{\x}_k;t_k)$ exists, as well as the time-derivative of the gradient $\nabla_{t \x} f(\hat{\x}_k;t_k)$, leading to the prediction model\footnote{The time derivative $\nabla_{t \x} f(\x;t_k)$ can be obtained via first-order backward finite difference if not available otherwise, see~\cite{Paper1, Bastianello2020}. }
\begin{eqnarray}\label{eq:model:un}
\!\!\!\varphi_k(\hat{\x}_{k+1|k}) = {\bf 0} \implies &\\	\hat{\x}_{k+1|k} = \hat{\x}_k - \nabla_{\x\x} ^{-1}f(\hat{\x}_k&\hspace*{-0.35cm};t_k)[\nabla_{\x} f(\x;t_{k}) + h\,\nabla_{t \x} f(\hat{\x}_k;t_k)].\hspace*{-.5cm}\nonumber
\end{eqnarray}
The prediction~\eqref{eq:model:un} represents a nonlinear discrete-time model to compute ${\x}_{k+1|k}$. Note that $\varphi_k({\x})$ \revi{can be interpreted as a specific choice for the gradient approximator $\mathbold{m}_k$ in~\cite{Jadbabaie2015}-- see discussion in the oracles paragraph}. 
Let us now consider a slightly more general setting than Problem~\eqref{eq:tip} as:
\begin{equation}\label{eq:fb}
\x^{\star}(t_{k}) = \argmin_{\x\in\reals^n}  f(\x;t_{k}) + g(\x) 
\end{equation}
where $g: \reals^n \to \reals\cup \{+\infty\}$ is a convex closed and proper function \rev{(e.g., $g(\x) = \|\x\|_1$)}. Problem~\eqref{eq:tip} is a special case of~\eqref{eq:fb}, when $g(\x)$ is the indicator function of the set $X$. Once again, we look at the first-order optimality conditions at time $t_{k}$, which can be framed as the generalized equation\rev{~\cite{Dontchev2009}}
\begin{equation}\label{eq:gen-equation}
	\nabla_{\x} f(\x;t_{k}) + \partial g(\x) \ni {\bf 0}.
\end{equation}
To predict how this first-order optimality condition changes in time (with information up to $t_k$), one can use a Taylor expansion around $(\hat{\x}_k; t_k)$, leading to the prediction model
\begin{equation}\label{eq:approx-gen-equation}
	\varphi_k(\hat{\x}_{k+1|k}) + \partial g(\hat{\x}_{k+1|k}) \ni {\bf 0}.
\end{equation}
Thus the prediction step requires the solution of this approximated generalized equation with initial condition $\hat{\x}_k$, which can be obtained, or approximated, cheaply with e.g., a few passes of a proximal gradient method~\cite{Eckstein1989} (cheaply since $\varphi_k$ is a quadratic function). The formulation~\eqref{eq:approx-gen-equation} represents the prediction model for the presented class of optimization problems~\eqref{eq:fb}, for a first-order Taylor expansion. Other prediction models exist for other classes of optimization problems~\cite{Paper1,Paper4}, for higher-order Taylor expansions~\cite{Dontchev2013}, and for more complex numerical integration methods~\cite{Jin2016, Liao2016, Guo2018, Qi2019, Qiu2019}. 

\smallskip

$\bullet$ Prediction based on parameter-estimation~\cite{Charles2016}. When the time dependence hides a parameter dependence, then models obtained via filtering are a viable alternative. Let $\b(t) \in \reals^l$ be a parameter, and let the function $f(\x;t) = f(\x;\b(t))$: e.g., the cost depends on the data stream $\b(t)$ representing for example the position of a robot to track. Then $\b(t)$ at time $t_{k+1}$ can be estimated via, for example, a Kalman filter based on the linear time-invariant model: 
\begin{equation}\label{eq:model}
\b(t_{k+1}) = \mathbold{\Gamma} \b(t_{k}) + \w_k, \quad  \y_k = \mathbold{\Phi} \b(t_{k}) + \mathbold{n}_k,
\end{equation}
for given matrices $\mathbold{\Gamma} \in\reals^{l\times l}$, $\mathbold{\Phi} \in \reals^{q\times l}$, observations $\y_k \in\reals^q$, and noise terms $\w_k \in \reals^l, \mathbold{n}_k \in \reals^q$. Then the prediction model requires the (approximate) solution of the problem
\begin{equation}
\hat{\x}_{k+1|k} \approx \argmin_{\x \in X} f(\x; \hat{\b}_{k+1}), 
\end{equation}  
with $\hat{\b}_{k+1}$ being the forecasted $\b(t_{k+1})$ based on the model~\eqref{eq:model} via e.g., a Kalman filter. Other models can be thought of based on non-linear models, more complicated forecasters, and even neural networks. 

%\smallskip
%$\bullet$ Other-prediction correction \as{Any other here?} 

\smallskip

\noindent {\bf Prediction-correction algorithms.} We have presented a few predictors for discrete-time time-varying optimization algorithms. No general result exists to encompass all the predictors. % but typically the results are based on how accurate is the prediction $\hat{\x}_{k+1|k}$. 
However, for a particular class of predictors (the one that employs first-order optimality conditions as prediction model) some general results can be derived. These methods are known as prediction-correction methods (since they predict how the optimization problem changes and then they correct for the errors in predictions based on the newly acquired cost~\cite{Paper1, Simonetto2017a}) and have roots in non-stationary optimization~\cite{Polyak1987,Popkov2005}, parametric programming~\cite{Dontchev2009,Zavala2010,Dontchev2013,Kungurtsev2017}, and continuation methods in numerical mathematics~\cite{Allgower1990}. %They also resemble evolutionary variational inequalities~\cite{Cojocaru2005,Nagurney2006} and path-following methods in interior point solvers~\cite{Nesterov2012}. 

Consider Problem~\eqref{eq:fb} for simplicity (although arguments are generalizable). Let $\mathcal{P}$ be a predictor method that approximates $\hat{\x}_{k+1|k}$ based on~\eqref{eq:approx-gen-equation}, in a Q-linear convergent fashion: one application of $\mathcal{P}$ acting on $\hat{\x}_k$ delivers a $\hat{\x}_{k+1|k}'$ for which
\begin{equation}\label{lin1}
\|\hat{\x}_{k+1|k}' - \hat{\x}_{k+1|k}\| \leq \varrho_1 \|\hat{\x}_k - \hat{\x}_{k+1|k}\|, \quad \varrho_1 \in (0,1).
\end{equation}
E.g., $\mathcal{P}$ could be a proximal gradient algorithm\rev{, in which case:
\begin{equation}\label{pred}
\hat{\x}_{k+1|k}' = \textrm{prox}_{\alpha g}\{\hat{\x}_k - \alpha \nabla_{\x} \varphi_k(\hat{\x}_k)\},
\end{equation}
where $\textrm{prox}_{\alpha g}\{\cdot\}$ is the proximal operator for function $g$ and step-size $\alpha$, } \revi{which could be applied one or multiple, say $P$, times for time step.} Let now $\mathcal{M}$, belonging to the same algorithm class of~\eqref{eq:linear}, be applied to the update (correction) step after function acquisition at $t_{k+1}$, for which one application on $\hat{\x}_{k+1|k}'$, delivers
\begin{equation}\label{lin2}
\|\hat{\x}_{k+1}-\x^{\star}(t_{k+1}) \| \leq \varrho_2 \|\hat{\x}_{k+1|k}'-\x^{\star}(t_{k+1}) \|, \quad \varrho_2 \in (0,1),
\end{equation}
\rev{for example another proximal gradient step as
\begin{equation}\label{corr}
\hat{\x}_{k+1} = \textrm{prox}_{\alpha g}\{\hat{\x}_{k+1|k}' - \alpha \nabla_{\x} f(\hat{\x}_{k+1|k}'; t_{k+1})\}.
\end{equation}
Then} the following result is in place.  

\begin{theorem}\label{th.2}\emph{(Informal)}
Consider the time-varying Problem~\eqref{eq:fb} and two methods $\mathcal{P}$ and $\mathcal{M}$ for which~\eqref{lin1}-\eqref{lin2} hold. Let the predictor $\mathcal{P}$ be applied $P$ times during the prediction step, and the corrector $\mathcal{M}$ be applied $C$ times. Consider Assumption~\ref{as:2} to hold and additionally, let $f(\x;t)$ be $L$-smooth (in addition to be $m$-strongly convex), with a well-defined Hessian $\nabla_{\x\x}f(\x;t)$. 
Then, there exists a minimal number of prediction and correction steps $P, C$ for which globally (i.e., starting from any initial condition)
$$
\limsup_{k \to \infty} \|\hat{\x}_{k} - \x^{\star}(t_{k})\| = \rev{\frac{\Delta_0}{m} \varrho_1^C \, O(h).}
$$
In addition, if we consider the assumption that higher-order derivatives of the cost function are bounded\revi{\footnote{\revi{Where induced Euclidean norms are considered for tensors.}}} as
$$
\max\{\|\nabla_{\x\x\x} f(\x; t)\|, \, \|\nabla_{t\x\x} f(\x; t)\|, \, \|\nabla_{tt\x} f(\x; t)\|\}\leq \rev{\Delta_1},
$$
uniformly in time and for all $\x \in \reals^n$, then locally (and for small $h$), there exists a minimal number of prediction and correction steps $P, C$ so that
$$
\limsup_{k \to \infty} \|\hat{\x}_{k} - \x^{\star}(t_{k})\| = \rev{\underbrace{O(\Delta_1 \varrho_1^C \,h^2)}_{\textrm{prediction gain}} + \underbrace{O(\Delta_0 \varrho_1^C \varrho_2^P \, h)}_{\textrm{approximation error}}.}
$$ 
\end{theorem}

\rev{
\begin{proof}
\emph(Sketch) The proof here proceeds as follows: we first bound the error coming from the prediction (e.g.,~\eqref{pred}), next bound the one from the correction (e.g.,~\eqref{corr}), and then combine them. For the error coming from the prediction, two errors must be considered, one coming from the model (due to the Taylor expansion error), the other coming from the $P$ prediction steps. When considering exact prediction ($P \to \infty$), the leading error is the Taylor expansion error (namely the error in~\eqref{eq:split}), which is $O(h)$ in general, and $O(h^2)$ when higher-order derivatives are bounded.   
\end{proof}
}

\smallskip

The results of Theorem~\ref{th.2} are fairly general and apply to different problem classes~\cite{Paper1, Paper4}. Theorem~\ref{th.2} indicates that tracking is not worse than correction-only methods in the worst case. If the function has some higher degree of smoothness, and we are interested in a local result, then a better ATE can be achieved, provided some (stricter) conditions on the number of prediction and correction steps are verified. The ATE is composed of two terms; one which is labeled as approximation error, which is due to the early termination of the prediction step (if $P \to \infty$ and prediction is exact, this term goes to $0$). The other, named prediction gain, is the gain coming from using a prediction step, which brings the error down to a $O(h^2)$ dependence on $h$. This depends on the first-order Taylor expansion employed; other methods can further reduce this to $O(h^4)$ or less~\cite{Dontchev2013, Jin2016, Liao2016, Guo2018, Qi2019, Qiu2019} (look again at Figure~\ref{fig:example}, \rev{where we have also employed a Taylor model up to degree $2$ for~\eqref{eq:split} to obtain an $O(h^3)$ error bound)}. 

The higher degree of smoothness required for the local results imposes boundedness of the tensor $\nabla_{\x\x\x}f(\x;t)$, which is a typical assumption for second-order algorithms (notice that the predictor requires second-order information, cf.~\eqref{eq:model:un}-\eqref{eq:approx-gen-equation} and its solution is comparable to solving a Newton step, which is \emph{locally} quadratically converging). Moreover, it bounds the variability of the Hessian of $f$ over time, which guarantees the possibility of performing more accurate predictions of the optimal trajectory.
Theorem~\ref{th.2} depicts a key result in prediction-correction methods: the prediction value is fully exploited with higher smoothness.% (since the predictor is based on Taylor expansion models). 

%%%%%%%%%%%%%%%%%%%%%%%%%%%%%%%%%%%%%%%%%%%%%%%%%%%
%%% SUBSECTION: CONTINUOUS-TIME ALGORITHMS
%%%%%%%%%%%%%%%%%%%%%%%%%%%%%%%%%%%%%%%%%%%%%%%%%%%

\subsection{Continuous-time algorithms} \label{sec_continuous}

%{\color{red} Is there a reason of using $\partial / \partial t$, instead of $\nabla_t$?}
%{\color{blue} {SP: Not really. } }

%{\color{red} any place for these references here? Continuous-time platforms have been discussed, e.g., in~\cite{Ye2015,Rahili2015a,Rahili2015,Gong2016,Ye2015,Xi2016a,Sun2017}.}

We consider now continuous-time prediction-correction algorithms which in general are appropriate in control and robotics applications\revi{\footnote{\revi{For these algorithms, time metrics like TR make less sense than in discrete-time setting. However, continuous-time algorithms are still interesting to investigate, both in theory: as continuous limits to discrete-time algorithms, and in practice: as good approximation of cases in which the sampling time is much smaller than other system characteristic times}.}}. The main component of these algorithms is the ability to track the minimizer by taking into account its evolution with time. {\color{black} {In continuous-time this scheme has been used in distributed time-varying convex optimization e.g.~\cite{Rahili2015,Gong2016,Xi2016a,Sun2017}.} } Since the objective function is $m$-strongly convex, the solution of the problem can be computed by solving the first-order optimality condition \eqref{eq:fermat}: for the implicit function theorem the time derivative of $\x^\star(t)$ is
  \begin{equation}\label{eqn_continuous_variation}
    \dot{\x^\star}(t) =  -\nabla^2_{\x\x} f(\x;t)^{-1}\nabla_{t\x}f(\x;t).
  \end{equation}
  In cases where the problem of interest is static, gradient descent and Newton's method can be used for instance to find trajectories such that $\lim_{t\to\infty} \x(t) = \x^\star$. Moreover, this convergence is exponential, meaning that there exist positive constants $C_1$ and $\alpha_1$ such that  $\left\|\x(t)-\x^\star\right\| \leq C_1 e^{-\alpha_1 t}$ (see e.g.,\cite[Definition 4.5]{khalil2002nonlinear} -- \revi{notice that exponential convergence is the continuous counterpart of the discrete-time $Q$-linear convergence}). To provide the same guarantees in the case of time-varying optimization we include the prediction term~\eqref{eqn_continuous_variation}, which incorporates changes in the optimizer
  \begin{equation}\label{eqn_dynamics}
\dot{\x}(t) = -\nabla^2_{\x\x} f(\x;t)^{-1}\left(\kappa\nabla_{\x}f(\x;t)+\nabla_{t\x}f(\x;t) \right),
    \end{equation}
  where $\kappa>0$ is \revi{referred to as ``gain of the controller'' in the literature, and~\eqref{eqn_dynamics} is referred to as ``the controller'', since it controls how the decision trajectory must change to reach the optimal solution trajectory.} %The previous differential equation defines a non-autonomous system for which exponential stability results exist. 
%  {\color{red} I am writing here some suggestions to streamline this, but they are not carved in stone. Suggested text: }
%  
%In particular, one can rely on standard results in the literature,  e.g.~\cite[Theorem 4.10]{khalil2002nonlinear}, to establish that $\e(t)=\x(t)-x^\star(t)$  converges exponentially to zero \cite{Zhao1998}, \cite[Proposition 1]{Fazlyab2018}.  
%%
%% THEOREM
%%
%\begin{theorem}
%Under the hypothesis of Theorem \ref{th.2}, $\x(t)$\textemdash the solution of the dynamical system \eqref{eqn_dynamics}\textemdash converges exponentially to $\x^\star(t)$, solution to \eqref{tvp}.
%
%    \begin{equation}
 %     \left\| \x(t)-\x^\star(t)\right\| \leq C_1e^{-Kt}
  %  \end{equation}
%  \end{theorem}
%  
%Maybe talk about Lyapunov here ? 
%
%
%  
%An interesting intuition is that by considering the error $\e(t) = \x(t)-\x^\star(t)$ one can show that the derivative of the Lyapunov function $V(\x;t) = \left\|\nabla_\x f(\x;t)\right\|^2/$ with respect to the time is zero. The latter is only possible since we are introducing the time variation of the optimal solution in our control scheme. If one does not take into account that term the tracking would be possible up to an asymptotic error that depends on the variation of the gradient with the time and the gain of the controller. 
%  
%  
%
This differential equation defines a non-autonomous dynamical system which converges exponentially to $\x^\star(t)$ \cite{Zhao1998}, \cite[Prop.~1]{Fazlyab2018}.  
  
%%
%% THEOREM
%%
\begin{theorem}\label{theo_tv}
Under the hypothesis of Theorem \ref{th.2}, $\x(t)$\textemdash the solution of the dynamical system \eqref{eqn_dynamics}\textemdash converges exponentially to $\x^\star(t)$, solution to \eqref{tvp}.
%
%    \begin{equation}
 %     \left\| \x(t)-\x^\star(t)\right\| \leq C_1e^{-Kt}
  %  \end{equation}
  \end{theorem}

\rev{
  \begin{proof}
    \emph(Sketch) The proof uses a Lyapunov argument. Define the error $\e(t) := \x(t)-\x^\star(t)$ and the function $V(\e;t) = \left\|\nabla_\x f(\e+\x^\star(t);t)\right\|^2/2$. Then the proof relies on establishing that $\dot{V}(\e;t) < 0$ for all $\e\neq 0$ and $\dot{V}(\mathbf{0};t)=0$ (see e.g.~\cite[Theorem 4.10]{khalil2002nonlinear}), and in particular: 
    %%
   % \begin{equation}
   %   \dot{V}(\e;t) = \nabla_\e V(\e;t)^\top\dot{\e} + \frac{\partial V(\e,t)}{\partial t}. 
   % \end{equation}
    %
    %Notice that the first term yields
    %
   % \begin{equation}
 %\nabla_\e V(\e,t)^\top\dot{\e} = \nabla_\x f(\e+\x^\star(t),t)^\top \nabla_{\x\x} f(\e+\x^\star(t),t)\dot{\e}.
 %   \end{equation}
    %
 %   On the other hand the second term yields
 %   \begin{align}
 %     \frac{\partial V(\e,t)}{\partial t}&= \nabla_\x f(\e+\x^\star(t),t)^\top \nabla_{\x\x} f(\e+\x^\star(t),t)\dot{\x}^\star(t)\nonumber \\
%      &+ \nabla_\x f(\e+\x^\star(t),t)^\top\frac{\partial \nabla_\x f(\e+\x^\star(t),t)}{\partial t}. 
%      \end{align}
    %
%    Combining the previous expressions and substituting $\e=\x-\x^\star(t)$ yields
%    \begin{equation}
%            \dot{V}(\e,t) = \nabla_\x f(\x,t)^\top \left(\nabla_{\x\x} f(\x,t)\dot{\x} +\frac{\partial \nabla_\x f(\x,t)}{\partial t} \right).
%    \end{equation}
%    Replacing $\dot{\x}$ by \eqref{eqn_dynamics} it follows that 
        \begin{equation*}
            \dot{V}(\e;t) = -\kappa \left\|\nabla_\x f(\e+\x^\star(t);t)\right\|^2 \leq 0.\vspace*{-15pt}
        \end{equation*}
%        Which is always negative for $\e\neq \mathbf{0}$ and equal to zero for $\e=\mathbf{0}$. Thus completing the proof of the result.
  \end{proof}
}

\smallskip

\rev{This result indicates that %In particular, by defining the error  $\e(t) = \x(t)-\x^\star(t)$ the function $V(\e,t) = \left\|\nabla_\x f(\e+\x^\star(t),t)\right\|^2/2$ is a Lyapunov function for the non-autonomous dynamical system $\dot{\e} = \dot{\x}-\dot{\x}^\star$.
the convergence is exponential to the optimal trajectory (ATE is \emph{zero}).} The latter is achieved by including the prediction in the controller, i.e., the time variation of the optimal solution. Without such predictor, tracking would be possible only up to an asymptotic error that depends on the variation of the gradient with the time and the gain of the controller. \rev{This is a clear benefit of structured algorithms. }
Notice that these results are the continuous time counter part of the results presented in Theorem \ref{th.2}. However, one of the advantages of working with continuous time flows is that it is also possible to establish asymptotic convergence to the solution of constrained optimization problems using interior point methods (see e.g., \cite[Chapter 11]{Boyd2004a}). Formally, let us define the following optimization problem
  \begin{subequations}\label{eqn_constrained}
  \begin{align}
      \x^\star(t) := &\argmin_{\x\in\mathbb{R}} && f(\x;t)\\
    &\mbox{subject to} && h_i(\x;t)\leq 0 \quad \forall i=1,\ldots,p.
    \end{align}
  \end{subequations}
  In \cite{Fazlyab2018}, inspired by interior point methods, the following barrier function is proposed
  \begin{equation}\label{eqn_barrier}
  \Phi(\x;t) = f(\x;t) -\frac{1}{c(t)}\sum_{i=1}^p \log\left(s(t)-h_i(\x;t)\right),
  \end{equation}
  where $c(t)$ is an increasing function such that $\lim_{t\to\infty}c(t)=\infty$ and $s(t)=s(0)e^{-\gamma t}$ for some $\gamma >0$. The intuition behind the barrier is that it approximates the indicator function as $t$ increases. This means that it takes the value $0$ when the constraint is satisfied and $+\infty$ in the opposite case. In that sense, when minimizing the unconstrained objective $\Phi(\x;t)$ constraint satisfaction is promoted. Notice that for the logarithm to be well defined we need $s(t)>h_i(\x;t)$ and thus the slack $s(t)$ is introduced just to guarantee that this is the case at all times $t\geq0$. In particular, it suffices to choose $s(0)\geq \max_{i=1,\ldots,p}\{h_i(\x(0),0)\}$ for this to be the case \cite[Theorem 1]{Fazlyab2018}. The previous intuition on how minimizing the function $\Phi(\x;t)$ defined in \eqref{eqn_barrier} resembles to solve \eqref{eqn_constrained} can be formally established. Let $\hat{\x}(t)$ be the minimizer of $\Phi(\x;t)$. Then it follows that $\lim_{t\to\infty}\left\|f(\hat{\x}(t);t)-f(\x^\star(t);t)\right\| =0$ \cite[Lemma 1]{Fazlyab2018}. This result, along with the idea that the barrier function promotes constraint satisfaction suggests that to solve \eqref{eqn_constrained} it suffices to compute the minimizer of the unconstrained barrier function $\Phi(\x;t)$ defined in \eqref{eqn_barrier}. %The latter can be done using the dynamical system discussed in \eqref{eqn_dynamics}. 
  This result is formalized in the following theorem.

  \begin{theorem}[ Theorem 1 \cite{Fazlyab2018}]\label{th.4}
    Consider the constrained optimization problem defined in \eqref{eqn_constrained} with $f(\x;t)$ $m$-strongly convex, $h_i(\x;t)$ for all $i=1,\ldots,p$ are convex functions and Slater's constraint qualifications hold: \rev{that is, there exists $\x^\dagger(t)$ such that for all $t\geq 0$ and for all $i=1,\ldots,p$ it holds that $h_i(\x^\dagger(t),t)<0$.} Let $\Phi(\x;t)$ be the barrier defined in \eqref{eqn_barrier} and let $\x(t)$ be the solution of the dynamical system:
    \begin{equation*}
      \dot{\x}(t) = -\nabla_{\x\x}\Phi(\x;t)^{-1}\left(\kappa\nabla_{\x}\Phi(\x;t) + \nabla_{t \x}\Phi(\x;t)\right).
    \end{equation*}
    Then it follows that  
    %
    %\begin{equation*}
    $
      \lim_{t\to\infty}\left\|\x(t)-\x^\star(t)\right\| = 0.
      $
   %   \end{equation*}
    \end{theorem}
\rev{
\begin{proof}
  The proof follows that of Theorem \ref{theo_tv} with $\e :=\x-\x^\star(t)$ and Lyapunov function $V(\e;t) = \left\| \nabla_\x \Phi(\e+\x^\star(t);t)\right\|^2/2$.
\end{proof}
}
\smallskip

Working in continuous time allows us to solve constrained problems using interior point methods, thus guaranteeing feasibility for all time\textemdash if the initial solution is feasible. This is especially appropriate for control systems where the constraints might represent physical constraints that need to be satisfied for the system to operate without failure. %In {\color{black}{Section VI-B}} we show how this barrier controller can be used to ensure safe navigation. 

% !TEX root = revised_final.tex

%%%%%%%%%%%%%%%%%%%%%%%%%%%%%%%%%%%%%%%%%%%%%%%%%%%%%%%%%%%%%%%%%%%%%%%%%%%%%%%%
%%%%%%%%%%%%%%%%%%%%%%%%%%%%%%%%%%%%%%%%%%%%%%%%%%%%%%%%%%%%%%%%%%%%%%%%%%%%%%%%
%%%%%%%%%%%%%%%%% APPLICATIONS
%%%%%%%%%%%%%%%%%%%%%%%%%%%%%%%%%%%%%%%%%%%%%%%%%%%%%%%%%%%%%%%%%%%%%%%%%%%%%%%%
%%%%%%%%%%%%%%%%%%%%%%%%%%%%%%%%%%%%%%%%%%%%%%%%%%%%%%%%%%%%%%%%%%%%%%%%%%%%%%%%

\section{Applications}\label{sec:applications} 

We highlight now application domains where structured and unstructured time-varying optimization methods have been or could be applied to. % solve relevant real-time information management, monitoring, control, and optimization tasks. 
%We also  provide pointers  to potential additional  applications for time-varying optimization tools and methods. %We focus on what is time-varying, what is the time scale and the typical computational complexity of the optimization problem. We also indicate whether structure optimization problems are or could be of use (i.e., a model, or oracle are easily available).  
%%
%It should be pointed out that the boundaries among the application areas  listed  below are difficult to delineate, given the increasingly cross-disciplinary nature of the research efforts on those domains and the natural overlaps  among different fields and communities  (with, e.g., communication research being important for control theory, as well as machine learning playing such an important role in signal processing and decision making). 
We proceed with a high-level (and by no means exhaustive) list of areas, presented in alphabetical order. Notice that, given the increasingly cross-disciplinary nature of the research efforts, clear boundaries are difficult to delineate.

\smallskip %%%%%%%%%%%%%%%%%%%%%%%%%%%%%%%%%%%%%%%%%%%%%%%%%%%%%%%%%%%%%%%%%%%% 

\noindent \textbf{Communications.} Problems such as congestion control, resource allocation, and power control have been of paramount importance in communication networks~\cite{Kelly1998, Low2002}. %, Srikant2004, Movsichoff2007}. 
Indeed, important questions arise when  channel capacities and non-controllable traffic flows are time-varying, with changes that are faster than  the solution time  of underlying optimization tasks, and even more so in the 5G era~\cite{Mehrabi2018} (for e.g., HD video streaming).  This setting can be tackled with time-varying optimization tools. %, with structured or unstructured online algorithms. % that can solve congestion control, resource allocation, and power control problems on the fly. 
For example in~\cite{Su2009}, a continuous-time structured algorithm with a first-order Taylor predictor model is proposed. %In the 5G era, with massive high-frequency data streaming (e.g., HD videos), over continent-size networks~\cite{Mehrabi2018,Park2019}, channel capacity time-variations and others will be more important and time-varying optimization will be even more prominent. 
The recent work~\cite{Bastianello2020} explored structured algorithms for intermittent time-varying service, a feature important in today's cloud computing. Finally time-variations are important when the communication graph is itself time-varying, see~\cite{Maros2019} and references therein.%, however this is slightly out of scope for the current overview and we refer to~\cite{Maros2019} for an up-to-date treatment and bibliography. %\as{Communication is not really my field: I guess somebody else should check this and expand. -- Geert? Yorgos? }

\smallskip %%%%%%%%%%%%%%%%%%%%%%%%%%%%%%%%%%%%%%%%%%%%%%%%%%%%%%%%%%%%%%%%%%%% 

\noindent \textbf{Control systems.} One popular tool in control systems is  model predictive control (MPC)~\cite{Mayne2000}. MPC is grounded on a strategy where an optimization problem is formulated to compute optimal states and commands for a dynamical system  over a given temporal window;  once a solution is identified, the command for the first time instant is implemented and the window is then shifted. The optimization problem changes over time, since it is parametrized over the state of a certain system, and it has to be re-solved every time. Recently, time-varying (and/or parameter-varying) algorithms for MPC have appeared for large-scale and embedded systems, e.g.,~\cite{Jerez2014,Hours2014,Gutjahr2016,paternain2019prediction}, which are a mix of continuous-time, discrete-time unstructured and structured algorithms. For example in~  \cite{Hours2014}, an unstructured algorithm (specifically a homotopy-based continuation method) is used to enhanced the tracking performance of the nonlinear MPC. In~\cite{paternain2019prediction} a predictor of the form \eqref{eq:model:un} is used to solve the optimization problem that arises from the receding horizon problem. Since the solution varies smoothly with the state of the system these methods are appropriate to achieve good tracking accuracy with low computational cost. %In particular, that of Inverting two Hessians. 
In \cite{liao2018semismooth}, these ideas are extended to problems with constraints by using a semi-smooth barrier function.

Other applications in control systems are the sequential training of neural networks for online system identification~\cite{Y.Zhao1993,Zhao1998,H.Myeong1994}, where predictors of the form  \eqref{eqn_dynamics} were proposed, as well as recent work at the intersection of online optimization and feedback control, where the output regulation problem is revisited by posing the problem of driving the output of a dynamical system to the optimal solution of a time-varying optimization problem~\cite{colombino2019online,zheng2019implicit}.

%Cross-fertilizing online optimization and feedback control, a recent line of works revisited the classical output regulation problem by posing the problem of driving the output of a dynamical system to the optimal solution of a time-varying
%optimization problem~\cite{colombino2019online,zheng2019implicit}.

\smallskip %%%%%%%%%%%%%%%%%%%%%%%%%%%%%%%%%%%%%%%%%%%%%%%%%%%%%%%%%%%%%%%%%%%% 

\noindent \textbf{Cyber-physical systems.} Cyber-physical systems (CPS)~\cite{Kim2012} are  engineered systems with tightly integrated  computing, communication, and control technologies. Because of major technological advances,  existing CPSs (power systems, transportation networks, and smart cities just to mention a few) are evolving towards societal-scale systems operating in highly dynamic environments, and with a massive number of interacting entities. It is then imperative to revisit information processing and optimization tools to enable optimal and  reliable decision-making on time scales that match the dynamics of the underlying physical systems. Due to space limitations, we focus here on  power systems and  transportation systems. 

A time-varying  problem for power systems can capture variations at a second level in non-controllable loads and available power from renewables~\cite{DhopleNoFuel}; it can also accommodate dynamic pricing schemes. Time-varying  problem formulations (and related online algorithms) can be utilized for tasks such as demand response, optimal power flow (OPF), and state estimation.  Adopting a time-varying optimization strategy, the power outputs of distributed energy resources (DERs) can be controlled at the second level to regulate voltages and currents within limits in the face of volatility of renewables and non-controllable loads, and to continually steer the network operation towards points that are optimal based on the formulated time-varying problem. Examples of works  include real-time algorithms for voltage control, optimal power flow, as well as DER management for aggregators; see for example~\cite{commelec1,DallAnese2016,Hauswirth2017,Tang2017,Liu2017,Liu2018} and pertinent references therein. For some applications such as the demand response and the OPF, online algorithms have been designed to leverage measurements of constraints (e.g., voltages  violations) in the algorithmic updates~\cite{DallAnese2016, Bernstein2019,tang2018feedback} to  relax the sensing  requirements. Real-time measurements were used in a state estimation framework in~\cite{song2019dynamic}. {\color{black} We develop these ideas with an example in Section~\ref{sec.pwr} .}% will illustrate an example of time-varying optimization problem and related online algorithms for real-time management of DERs in power distribution grids. }

In the context of  transportation systems, fast time-variations may arise  from different factors (and at appropriate time-scales), such as variations in the traffic, pedestrians crossing the roads, car accidents, sport events; these factors may lead to  time-dependent routing and traffic light control algorithms~\cite{Gendreau2015}. Motivated by the recent widespread use of ridesharing and mobility-on-demand services~\cite{Alonso-Mora2017}, spatio-temporal-variations naturally emerge from customer pick-up and drop-off requests as well as fleet locations. As representative works in context,~\cite{Alonso-Mora2017} and~\cite{Simonetto2019} discussed unstructured algorithms to achieve long-term (``asymptotical'') good tracking, while sacrificing short term optimality. In~\cite{Alonso-Mora2017a}, an online algorithm based on a structured problem formulation is presented, where the prediction is based on historical data and machine learning forecasting. An unstructured algorithm is also presented in~\cite{Eser2018}, to find optimal meeting points.

\smallskip %%%%%%%%%%%%%%%%%%%%%%%%%%%%%%%%%%%%%%%%%%%%%%%%%%%%%%%%%%%%%%%%%%%% 

\noindent \textbf{Machine learning and signal processing.} As a representative problem spanning the broad fields of machine learning and signal processing,  we focus on the reconstruction of sparse signals via $\ell_1$-regularization where we are interested in recovering a sparse signal given some observations, e.g., extract ``sparse'' features in images~\cite{Beck2009}. The time-varying nature of this problem arises when we want for \rev{instance} to extract features in videos. Works that explore dynamic $\ell_1$ reconstruction are, e.g.,~\cite{Balavoine2015,Yang2015,Vaswani2015}. In~\cite{Asif2014, Charles2016}, two algorithms are presented, one unstructured using homotopy and one structured building a model based on methods akin to Kalman filters. In~\cite{Fosson2018}, unstructured methods for the elastic net are discussed. 

Other applications in machine learning and signal processing, where a number of (mainly) unstructured algorithms have been proposed, include contemporary approaches for sparse, kernel-based, robust, linear regression, zeroth-order methods, and learning problems over networks. Additional lines of work include  dynamic classification under concept drift~\cite{Das2016}, dynamic beamforming~\cite{Maros2017}, and other dynamic signal processing tasks, such as maximum a posteriori estimation~\cite{Ling2013,Jakubiec2013}.

\smallskip %%%%%%%%%%%%%%%%%%%%%%%%%%%%%%%%%%%%%%%%%%%%%%%%%%%%%%%%%%%%%%%%%%%% 

\noindent \textbf{Medical engineering.} Medical engineering is a growing research field in many contexts. Here we focus on the new possibilities offered by new and fast imaging modalities under magnetic resonance (see~\cite{Clogenson2016,Rueckert2019} for a broader context). Once confined to static images (due to the high computational load), MRI is now transitioning to fast imaging and possibly high definition video streaming, which could be of invaluable help to clinicians and researchers alike, not to mention patients, especially children. In the series of work~\cite{Zhang2010b, Lingala2016}, the authors describe an unstructured algorithm to solve a time-varying subsampled nonlinear regularized inverse problem. The algorithm allows the clinicians to visualize blood flow, cardiac features, and swallowing, among many other things. %Others in case we have space: {Ntsinjana2011, Voit2013, Olthoff2014, Zhang2014}

\smallskip %%%%%%%%%%%%%%%%%%%%%%%%%%%%%%%%%%%%%%%%%%%%%%%%%%%%%%%%%%%%%%%%%%%% 

\noindent \textbf{Optimization and mathematical programming. } Time-varying optimization has been studied for applications within mathematical programming, e.g., in the context of parametric-programming~\cite{Guddat1990,Zavala2010,Dinh2012,Dontchev2013,Kungurtsev2017} where a wealth of structured and unstructured algorithms are presented. Time-varying optimization has its roots in continuation methods in numerical mathematics~\cite{Allgower1990} and it resembles path-following methods~\cite{Nesterov2012}, so advances in either fields are intertwined. 

Another  application in mathematical programming where time-varying optimization could be (and has been) applied is the field of evolutionary variational inequalities. Variational inequalities~\cite{Kinderlehrer1980} can be framed as optimization problems, while evolutionary ones can be framed as time-varying optimization problems. In~\cite{Daniele2004,Cojocaru2005,Nagurney2006}, the authors discuss plenty of interesting applications in socio-economical sciences (human migration studies, economics, time-dependent equilibria in games, etc.), proposing mainly unstructured approaches. 

\smallskip %%%%%%%%%%%%%%%%%%%%%%%%%%%%%%%%%%%%%%%%%%%%%%%%%%%%%%%%%%%%%%%%%%%% 

\noindent \textbf{Process Engineering.} In chemical and process engineering, the body of work~\cite{Chachuat2009, Jaeschke2011,Graciano2015} focuses on real-time optimization for chemical and industrial processes. The optimization problem is not time-varying per se, but it becomes time-varying because the constraints (i.e., the industrial process) are learned online and adapted. Several real-time optimization algorithms are proposed, mainly unstructured. 

\smallskip %%%%%%%%%%%%%%%%%%%%%%%%%%%%%%%%%%%%%%%%%%%%%%%%%%%%%%%%%%%%%%%%%%%% 

\noindent \textbf{Robotics.} Time varying optimization problems\textemdash or problems that depend on a time-varying parameter\textemdash appear often in the context of robotic systems. In the context of safe navigation \cite{arslan2016exact,arslan2019sensor} consider the problem of using power diagrams to define a local safe space, which depends on the position of the agent itself. The control law used to navigate is such that it aims to track the projection of the goal on the local safe-space. Even in cases, where the goal is static, a time-varying optimization problem needs to be solved due to the modification of the local free space. In \cite{Fazlyab2018} the approach described in Section \ref{sec_continuous} is used to compute said solutions. We develop these ideas more in ~\ref{sec.robots}.

For networks of mobile robots \cite{Zavlanos2013}, the ``communication integrity'' is guaranteed by solving a time-varying optimization problem. Specifically, since an unstructured algorithm is used, an asymptotic tracking error that results in small constraint violation and sub-optimality is achieved. %Since the rates depend on the position of the agents and these change over time to perform the required task, the solution of a sequence of time-dependent\textemdash or position dependent to be more precise\textemdash optimization problems needs to be computed. Since the problem is constrained and the Lagrangian associated to the problem can be maximized in closed form, one can solve the problem in the dual domain by running gradient descent. Specifically, if a predictor is not used as it is the case in \cite{Zavlanos2013} an asymptotic tracking error that results in small constraint violation and sub-optimality is achieved. 

Another interesting application is that of robotic manipulators~\cite{Miao2015, Liao2016, Li2018}, obtained via zeroing neural dynamics (ZND) \cite{zhang2011zhang, zhang2015taylor, zhang2019zeroing}, based on a prediction step similar to \eqref{eqn_dynamics}. %In \cite{Miao2015, Liao2016, Li2018}, the online motion planning problem is cast as a time-varying quadratic objective linearly constrained optimization problem. Since the computation of these solutions needs to be done in the control interval, methods without predictors might take too-long to converge. Zeroing neural dynamics (ZND) \cite{zhang2011zhang} are based on a prediction step similar to \eqref{eqn_dynamics} and enjoy the advantage of allowing parallel computing. Their drawback is that, the proposed dynamics are continuous and in terms of practical implementations, discrete-time models are the basic considerations for common digital equipment. The main challenge in discretizing these dynamics is that traditional one-step ahead approaches\textemdash to make the implementation fully online \textemdash cannot be applied since they are unable to stabilize the resultant models \cite{zhang2015taylor}. The article \cite{zhang2019zeroing} provides a good introduction to other algorithms that allow to discretize the ZND.

%\smallskip %%%%%%%%%%%%%%%%%%%%%%%%%%%%%%%%%%%%%%%%%%%%%%%%%%%%%%%%%%%%%%%%%%%% 

%\noindent \textbf{Signal processing.} 

% !TEX root = revised_final.tex

%%%%%%%%%%%%%%%%%%%%%%%%%%%%%%%%%%%%%%%%%%%%%%%%%%%%%%%%%%%%%%%%%%%%%%%%%%%%%%%%
%%%%%%%%%%%%%%%%%%%%%%%%%%%%%%%%%%%%%%%%%%%%%%%%%%%%%%%%%%%%%%%%%%%%%%%%%%%%%%%%
%%%%%%%%%%%%%%%%% DEEP DIVES
%%%%%%%%%%%%%%%%%%%%%%%%%%%%%%%%%%%%%%%%%%%%%%%%%%%%%%%%%%%%%%%%%%%%%%%%%%%%%%%%
%%%%%%%%%%%%%%%%%%%%%%%%%%%%%%%%%%%%%%%%%%%%%%%%%%%%%%%%%%%%%%%%%%%%%%%%%%%%%%%%

\section{Two applications: deep dive}
\label{sec:deepdiveapplications}

\subsection{Example in Power Grids}\label{sec.pwr}

Consider a power distribution grid serving residential houses or commercial facilities, featuring  $N$ controllable DERs. The vector  $\x_i \in X_i \subset \mathbb{R}^2$ collects the active and reactive power outputs of the $i$th DER, and $X_i$ models hardware constraints. A prototypical time-varying optimization problem for real-time management of DERs is:
\begin{equation}
\label{eq:tip_sg}
\x^{\star}(t_k) := \argmin_{\{ \x_i \in X_i \}_{i = 1}^N} \, \sum_{i = 1}^N f_i(\x_i; t_k) + f_{N+1}(\x; t_k)
\end{equation}
where $f_i(\x_i; t_k)$ is a cost function associated with the $i$th DER, and $f_{N+1}(\x; t_k)$ is a time-varying cost associated with the power network operator. Elaborating on the latter, suppose for example  that a linearized model for the power flow equations is utilized to capture the variations on some electrical quantities $\y \in \mathbb{R}^m$ (e.g., voltages and power flows on lines) induced by $\x$; that is, $\y(t_k) = \A_x \x + \A_w \w(t_k)$, where $\w(t_k)$ is a vector collecting the powers of non-controllable devices and $\A_x, \A_w$ are sensitivity matrices that are built based on the network topology and the line impedances~\cite{Liu2018,DallAnese2016}. A possible choice for the function $f_{N+1}(\x; t_k)$ for the network operator can then  be $f_{N+1}(\x; t_k) = \frac{\gamma}{2}\|\y^{\mathrm{ref}}(t_k) -  \A_x \x + \A_w \w(t_k) \|^2$, where $\y^{\mathrm{ref}}(t_k)$ it a time-varying reference point for the electrical quantities included in $\y$, and $\gamma > 0$  is a design parameter that influences the ability to track the time-varying reference signal  $\y^{\mathrm{ref}}(t_k)$. Various models for $f_i(\x_i; t_k) $ can be adopted, based on specific problem settings; for example, $f_i(\x_i; t_k) = \|\x_i - \x_i^{\mathrm{ref}}(t_k) \|^2$ can minimize the deviation from a desirable setpoint for the $i$th DER (that can be computed based on a slower time-scale dispatch problem); in case of photovoltaic systems, $\x_i^{\mathrm{ref}}(t_k) $ could be set to $\x_i^{\mathrm{ref}}(t_k) = [P^{\mathrm{av}}(t_k), 0]^T$, with $P^{\mathrm{av}}(t_k)$ the maximum power available, to minimize the power curtailed. Alternatively, set  $f_i(\x_i; t_k)$ to a time-varying incentive  $-\!\bpi_i^T (t_k) \x_i$ to maximize the profit of the $i$th DER in providing services to the grid.%  $\pi_i (t_k)$ is a time-varying incentive. 

{ \color{black} 
With reference to Figure~\ref{fig:depict}, in this application data streams include the parameters of the time-varying function $f_i(\x_i; t_k)$ (e.g., the power setpoints $\{\x_i^{\mathrm{ref}}(t_k)\}$ or the incentive signals $\{\bpi_i (t_k)\}$), the function  $f_{N+1}(\x; t_k)$ (where set points $\y^{\mathrm{ref}}(t_k)$ can rapidly change to provide various services to the grid), as well as the powers $\w(t_k)$ consumed by the non-controllable devices. The algorithm produces decisions on setpoints for the  active and reactive power outputs  $\x_i(t_k)$ of the DERs, which are commanded to the devices. Finally, ``feedback'' can come in the form of measurements of the actual power outputs $\x_i(t_k)$~\cite{commelec1}, as well as other electrical quantities~\cite{DallAnese2016,tang2018feedback}.
}

\begin{figure}
\centering
\includegraphics[width= 0.425\textwidth, trim=0cm 0cm 0 1.25cm, clip=on]{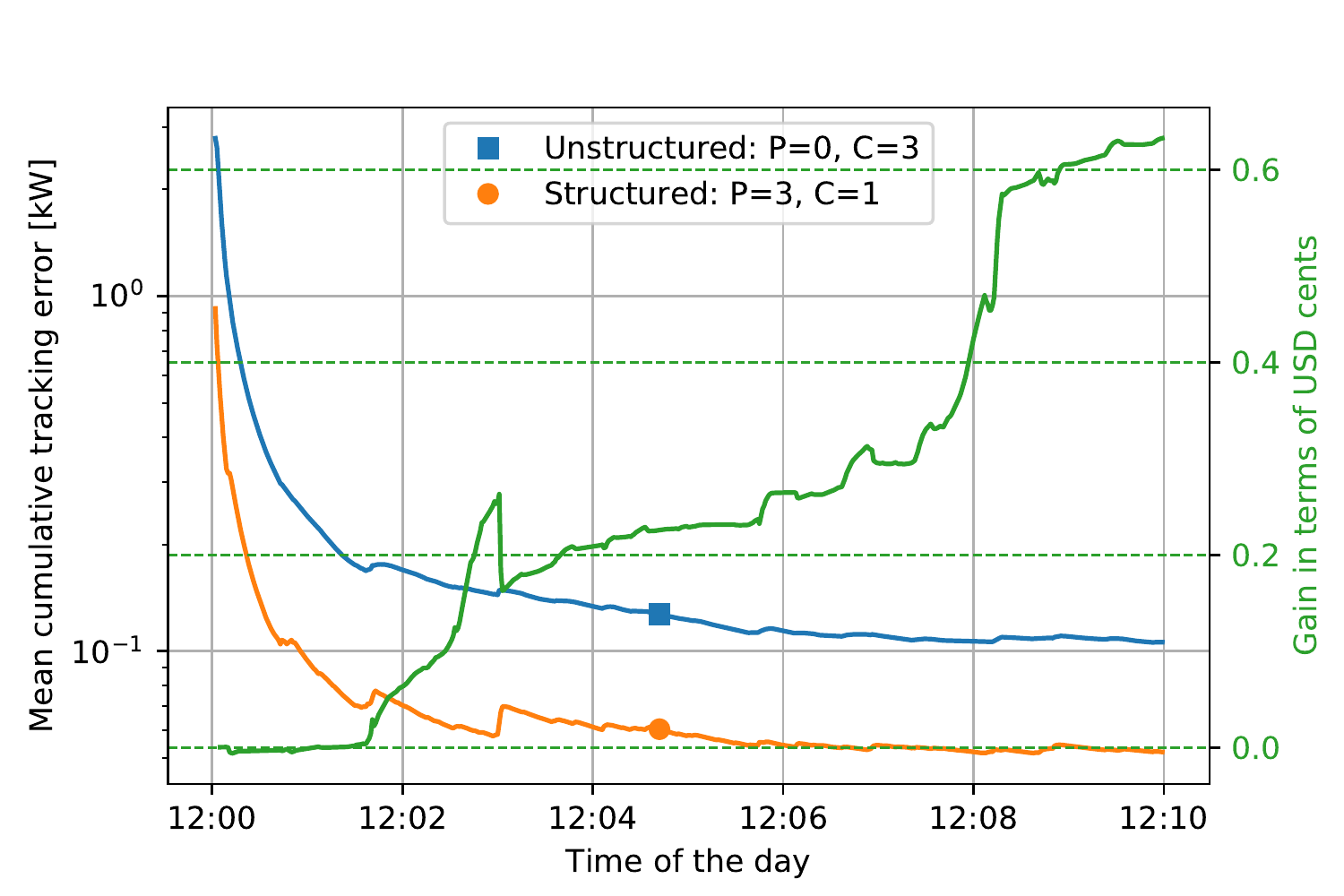}
\caption{\small {Mean cumulative tracking error $(1/T) \sum_{k = 1}^T \|\x(t_k) - \x^{\star}(t_k)\|$ vs. time of the day for a choice of structured ($P=3$, $C=1$) and unstructured ($C=3$) algorithms, having the same computational time. In green, we report the hypothetical gain in terms of less utilized power at the average cost of 12 USD cent per kWh.}}
\label{fig:example-smartgrid}
\end{figure}

As an illustrative example, we consider the case where $N = 500$ DERs are controlled in a distribution feeder; the set $X_i$ is build so that the ranges of active and reactive powers are $[-50, 50]$ kW and $[-50, 50]$ kVAr,  and $f_i(\x_i; t_k)$ is set to $f_i(\x_i; t_k) = \frac{1}{2} \|\x_i\|^2$ for all DERs. This setting is representative of a case where energy storage resources are utilized to provide services. We consider the case where $\y$ is a scalar and represents
the net power consumed by a distribution network; in this case, $\y^{\mathrm{ref}}(t_k)$ can model automatic generation control (AGC) signals or flexible ramping signals. The matrices $\A_x, \A_w$  are built as in~\cite{DallAnese2016}. We use the real data provided in~\cite{DallAnese2016} to generate the vectors $\w(t_k)$ with a granularity of one second. The parameters are $m = 1$ and $L = 21$, and $\gamma = 2$; the step size is $\rev{\alpha} = 1/(10 L)$. We keep the computational time fixed in our comparison between the \rev{unstructured} running projected gradient and \rev{the structured} prediction-correction algorithm; in particular, we consider the cases $P= 0, C = 3$ and $P = 3, C = 1$ (see Theorem~\ref{th.2}).  

\rev{
To outline the steps of the prediction-correction algorithm, recall that $\hat{\x}_k$ denotes the iterate of the algorithm at time $t_k$ [cf. Thm.~\ref{th.1}], and let $f(\x;t_k) :=  \sum_{i = 1}^N f_i(\x_i; t_k) + f_{N+1}(\x; t_k)$ and $X = X_1 \times X_2 \times \ldots  \times X_N$ for brevity. A prediction $\hat{\x}_{k|k-1}$ is obtained by running $P$ prediction steps $p = 0, \ldots, P-1$:
\begin{multline}
\hat{\x}^{p+1}  = \textrm{proj}_{X} \left\{ \hat{\x}^{p} - \alpha \left( \nabla_{\x\x}f(\hat{\x}_k;t_k)(\hat{\x}^{p} - \hat{\x}_{k-1}) \right. \right.   \\
 \left. \left. +  h \nabla_{t \x}  f(\hat{\x}_{k-1};t_k) + \nabla_{\x} f(\hat{\x}_{k-1};t_k) \right) \right\},
\label{eq:sg-pred}
\end{multline}
and by setting $\hat{\x}_{k|k-1} = \hat{\x}^{P}$. Starting now from $\bar{\x}^{0} =  \hat{\x}_{k|k-1}$, the correction phase involves the following $C$ steps: 
\begin{align}
\bar{\x}^{c+1} & = \textrm{proj}_{X} \left\{ \hat{\x}^{p} - \alpha \left( \nabla_{\x} f(\bar{\x}^c ;t_k) \right) \right\}
\label{eq:sg-corr}
\end{align}
for $c = 0,1, \ldots, C-1$. The  iterate $\hat{\x}_{k}$ is then  $\hat{\x}_{k} = \bar{\x}^{C}$.  Notice that if $P = 0$, one recovers the unstructured running projected gradient method; see also~\eqref{prj.gr}. In the simulations, the time derivative $ \nabla_{t \x}  f(\hat{\x}_k;t_k)$ in \eqref{eq:sg-pred} is substituted by an approximate version (see, for example,~\cite{Paper1} and~\cite{Paper2}).
}

To assess the performance of the prediction-correction algorithm,  Figure~\ref{fig:example-smartgrid} depicts the mean cumulative tracking error $(1/T) \sum_{k = 1}^T \|\x(t_k) - \x^{\star}(t_k)\|$. It can be seen that by leveraging the temporal structure of the problem, the prediction-correction algorithm offers improved performance. \rev{We can now evaluate the performance metrics presented in Section~\ref{sec:tvo}. We compute the ATE as the mean error in the last $20$~s of the simulation, yielding an ATE of $\sim50$~W for the unstructured method, and an ATE of $\sim80$~W for the structured method. Since the computational time of both methods is the same, it follows that SG $= 1.6$. The CR can be empirically evaluated by the time it takes to enter the ATE ball as approximately $1$~minute for both methods. On the other hand, the TR is hardware-dependent, since the denominator of the TR depends on the computational capabilities of the microcontrollers embedded in the DERs, where algorithms are implemented. } %\rev{Finally, we report on the hypothetical gain in terms of less  energy savings  at the average cost of 12 USD cent per kWh, further observing the importance of using structured methods.}

\subsection{Example in Robotics}\label{sec.robots}

%%%%%%%%%%%%%%%%%%%%%%%%%%%%%%%%%%%%%%%%%%%%%%%%%%%%%%%%%%%%%%%%%%%%%%%%%%%%%%%%%% F I G U R E %%%%%%%%%%%%%%%%%%%%%%%%%%%%%%%%%%%%%%%%%%%%%%%%%%%%%%%%%%%%%%%%%%%%%%%%%%%%%%%%%%%%%%%%%%%%%%%%%%%%%%%%%%%
\begin{figure}
\vspace{-20pt}
\includegraphics[width=0.52\textwidth]{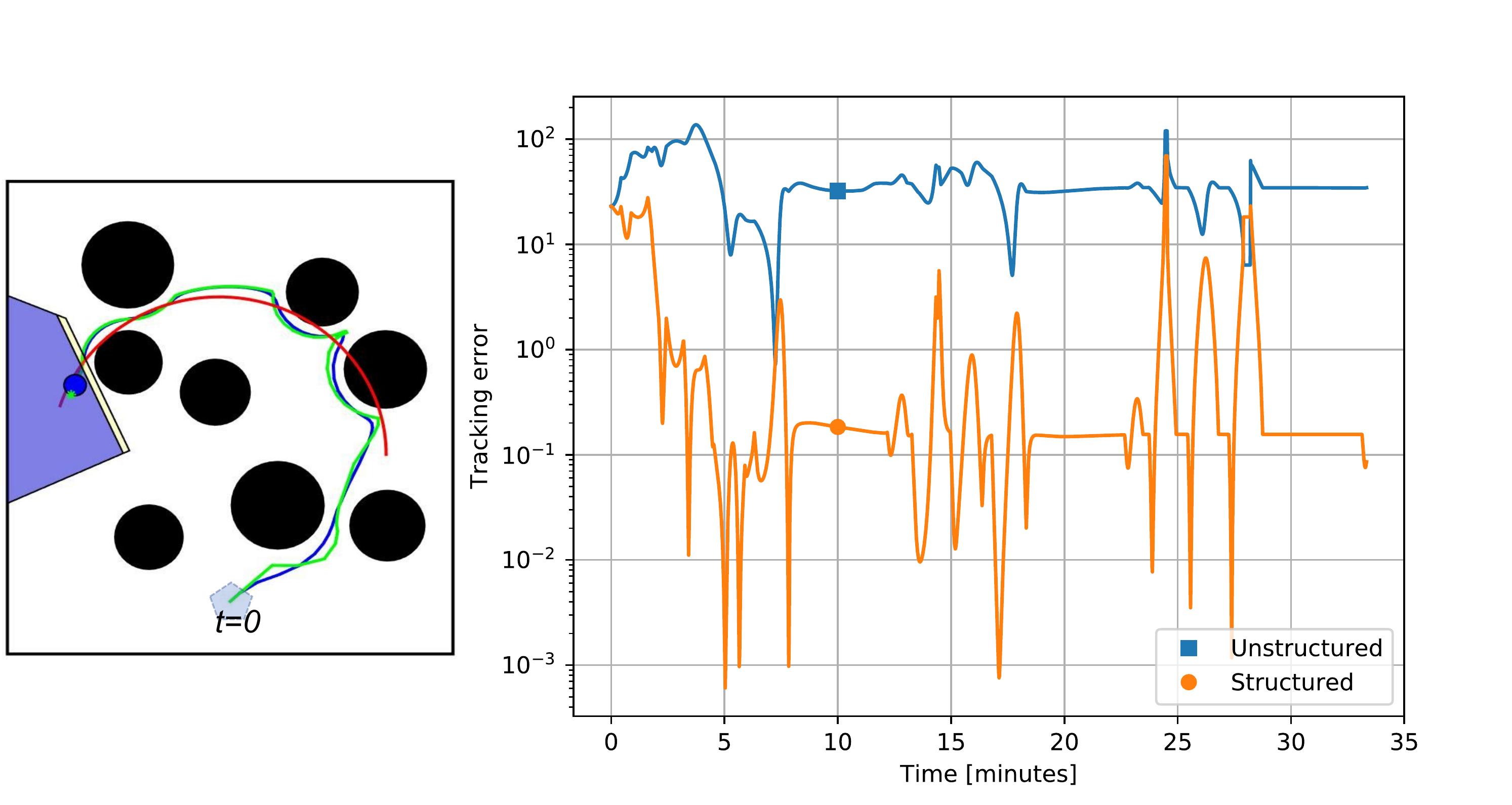}
\caption{{\small Left: The red circle represents the desired  $\x_{\textsf{d}}(t)$. The green and blue lines represent, respectively, the trajectories of the estimates of the projected goal $\widehat{\x}(t)$ and the trajectories of the robot $\x_{\textsf{r}}(t)$ for the structured algorithm. Right: Tracking error $\|\x(t_k) - \x^{\star}(t_k)\|$ vs. time for a choice of structured and unstructured algorithm.}}
\label{fig_dynamic_goal_trajectory}
\end{figure}
%%%%%%%%%%%%%%%%%%%%%%%%%%%%%%%%%%%%%%%%%%%%%%%%%%%%%%%%%%%%%%%%%%%%%%%%%%%%%%%%%% F I G U R E %%%%%%%%%%%%%%%%%%%%%%%%%%%%%%%%%%%%%%%%%%%%%%%%%%%%%%%%%%%%%%%%%%%%%%%%%%%%%%%%%%%%%%%%%%%%%%%%%%%%%%%%%%%

Consider a navigation setup of driving a disk-shaped robot of radius $r>0$, whose position is denoted by $\x_{\textsf{r}}(t)$, to a desired configuration $\x_{\textsf{d}}(t)$, while avoiding collisions with obstacles in the environment. Here, we deal with a closed and convex workspace $\mathcal{W} \subset \mathbb{R}^n$ of possible configurations that the robot can take. Assume that the workspace is populated with $m$ non-intersecting spherical obstacles, where the center and radius of the $i$-th obstacle are denoted by $\x_i \in \mathcal W$ and $r_i >0$, respectively. \rev{In general, this navigation problem is nonconvex due to the presence of obstacles, however one can convexify it by looking at the \emph{collision-free convex local workspace} around $\x_{\textsf{r}}$, defined as~\cite{arslan2016exact} }
%We start by defining the \emph{power distance} between a point $\x$ and a disk $B(\x_i,r_i )$ as $\mathcal{P}(\x,B(\x_i,r_i ))=\left\Vert\x-\x_i\right\Vert_2^2-r_i^2$, \cite{aurenhammer1987power}. Following the development in {\color{black}{[REF]}} denote the \emph{local workspace} around $\x_{\textsf{r}}$ as
%\[
%\mathcal L \mathcal W(\x_{\textsf{r}})=\left\{\x \in \mathcal W \colon \mathcal{P}(\x,B(\x_{\textsf{r}},r))\leq \mathcal{P}(\x,B(\x_i,r_i )),\, \forall i \right\},
%\]
%i.e., the set of points in $\mathcal{W}$ that are closer (in power distance) to the robot than to any of the obstacles.
%The local workspace defines a polytope whose boundaries are hyperplanes, such as the polygon marked with a thick light blue line in Fig. \ref{fig_static_goal_trajectory} (see \cite[Eq. (6)]{arslan_kod_ICRA2016B} for an explicit expression of these hyperplanes).
%Furthermore, the
%\emph{collision-free local workspace} around $\x_{\textsf{r}}$ as %\cite{arslan_kod_ICRA2016B}:
\begin{equation}
\mathcal L \mathcal F (\x_{\textsf{r}}) \!= \!\left\{\! \x\in \mathcal W \colon (\x_i\!-\!\x_{\textsf{r}})^\top \x - b_i(\x_{\textsf{r}}) \leq 0, \, i=1\ldots m\!\right\} \nonumber
\end{equation}
\rev{ %with $\a_i(\x_{\textsf{r}}) = \x_i-\x_{\textsf{r}}$, and 
%\begin{equation}
%\b_i(\x_{\textsf{r}}) = (\x_i-\x_{\textsf{r}})^\top\left(\theta_i\x_i+(1-\theta_i)\x_{\textsf{r}}+r\frac{\x_{\textsf{r}}-\x_i}{\left\|\x_{\textsf{r}}-\x_i\right\|}\right),
%\end{equation}
%
% where in the previous equation $\theta_i$ are scalars dependent on $\x_{\textsf{r}}$ defined as
 %
% \begin{equation}
% \theta_i = \frac{1}{2}-\frac{r_i^2-r^2}{\left\|\x_i-\x_{\textsf{r}}\right\|^2}.
% \end{equation}
where $b_i(\x_{\textsf{r}})$ are pertinent scalars computed depending on robot and obstacles positions (cf.~\cite{arslan2016exact})
}
\rev{The collision-free local workspace describes a local neighborhood of the robot that is guaranteed to be free of obstacles. Each obstacle introduces a linear bound, and thus the local free space is convex, and yields a polygon as the blue colored one in  Figure~\ref{fig_dynamic_goal_trajectory}} (cf. \cite[Eq. (6)]{arslan2016exact}). \rev{The position of the target $\x_{\textsf{d}}(t)$, the location of the robot $\x_{\textsf{r}}(t)$, and the local free space $\mathcal L \mathcal F (\x_{\textsf{r}})$ correspond to the data stream of Figure \ref{fig:depict}.}
Supposing that the robot follows the integrator dynamics $\dot{\x}_{\textsf{r}}= \mathbf{u}(\x_{\textsf{r}})$, the controller proposed in \cite{arslan2016exact} is given by 
%
%\begin{equation}\label{eqn_omur_controller}
$\dot{\x}_{\textsf{r}}(t) = -G_{\textrm{c}} ( \x_{\textsf{r}} - \x^\star(t)), 
$%\end{equation}
where $G_{\textrm{c}}>0$ and $\x^\star(t)$ is the orthogonal projection of the desired configuration $\x_{\textsf{d}}(t)$ onto the collision-free local workspace $\mathcal L \mathcal F(\x_{\textsf{r}})$. \rev{Since the local workspace is collision free, so is the direction $\x_{\textsf{r}}-\x^\star(t)$, and thus the control law is guaranteed to avoid the obstacles. This controller also guarantees that the robot converges to $\x_{\textsf{d}}(t)$ \cite{arslan2016exact}.} %
%\begin{figure}
%\hspace{-30pt}
%\includegraphics[width=0.60\textwidth]{./Figures/moving_target.jpeg}
%\caption{}
%\label{fig_ate_robot}
%\end{figure}
%
%%%%%%%%%%%%%%%%%%%%%%%%%%%%%%%%%%%%%%%%%%%%%%%%%%%%%%%%%%%%%%%%%%%%%%%%%%%%%%%%%% F I G U R E %%%%%%%%%%%%%%%%%%%%%%%%%%%%%%%%%%%%%%%%%%%%%%%%%%%%%%%%%%%%%%%%%%%%%%%%%%%%%%%%%%%%%%%%%%%%%%%%%%%%%%%%%%%
%\begin{figure}
%\centering
%\includegraphics[width=0.20\textwidth]{./Figures/trajectories3}
%\caption{{\small The red circle represents the desired configuration $\x_{\textsf{d}}$. The green and blue lines represent, respectively, the trajectories of the estimates of the projected goal $\widehat{\x}(t)$ and the trajectories of the robot $\x_{\textsf{r}}(t)$ for 4 different initial conditions.}}
%\label{fig_static_goal_trajectory}
%\end{figure}
%
%\rev{The proposed control law corresponds to the decision stream and feedback of Figure \ref{fig:depict}. 
\rev{It requires computating the projection of $\x_{\textsf{d}}(t)$ onto $\mathcal L \mathcal F(\x_{\textsf{r}})$ by \emph{solving the time-varying convex problem}} 
\begin{equation}\label{eqn_projected_goal}
\x^\star(t):= \argmin_{\x \in \mathcal L \mathcal F (\x_{\textsf{r}}) \subseteq \mathbb{R}^n} \frac{1}{2} \| \x-\x_{\textsf{d}}(t)\|^2. % \\
%&\mbox{s.t.} \quad \mathbf{a}_i(\x_{\textsf{r}})^\top\x - b_i(\x_{\textsf{r}}) \leq 0, \quad i =1\ldots m.
\end{equation}
By using the barrier function defined in \eqref{eqn_barrier} and the dynamics in Theorem \ref{th.4}, one can compute $\hat{\x}(t)$, an estimate  of $\x^\star(t)$, and apply the control law %\eqref{eqn_omur_controller}, 
%with the difference that we use an estimate of the projected goal instead of the projected goal itself, i.e., we consider the closed loop dynamics
%
%\begin{equation}\label{eqn_control_law}
$\dot{\x}_{\textsf{r}}(t) = -G_{\textrm{c}}(\x_{\textsf{r}}-\widehat{\x}(t))$.
%\end{equation}
%
%
\rev{
The barrier function in \eqref{eqn_barrier} for this application takes the form 
\begin{equation*}
\Phi(\x,\x_{\textsf{r}};t) = \frac{1}{2}\left\|\x-\x_{\textsf{d}}(t)\right\|^2-\frac{1}{c(t)}\sum_{i=1}^m \log(b_i(\x_{\textsf{r}})-\a_i(\x_{\textsf{r}})^\top\x)
\end{equation*}
with $\a_i \!=\! \x_i \!-\! \x_{\textsf{r}}$. Then estimate $\hat{\x}(t)$ is the solution of the following dynamical system with initial condition $\hat{\x}(0) = \x_{\textsf{r}}(0)$
\begin{equation*}
\dot{\hat{\x}}(t) \!=\! -\!\nabla_{\x \x}\Phi(\x,\x_{\textsf{r}};t)\!^{-1}\!\!\left(\kappa\nabla_{\x}\Phi(\x,\x_{\textsf{r}};t)\!+\!\nabla_{\x t}\Phi(\x,\x_{\textsf{r}};t)\!\right),
\end{equation*}
where $c(t)= 1 e^{0.1 t}$ and $\kappa =0.1$.
}
\rev{
To evaluate the performance of the proposed controller and optimizer, we consider a workspace $\mathcal W=[-20,20]\times[-20,25]$ containing  $8$ circular obstacles (black circles in Figure~\ref{fig_dynamic_goal_trajectory}-left). Figure~\ref{fig_dynamic_goal_trajectory}-left, also depicts the trajectories followed by a disc-shaped robot of radius equal to one (blue circle) where $G_{\textrm{c}}=2$. The red line represents the trajectory of $\x_{\textsf{d}}(t)$ and the green and blue lines represent, respectively, the trajectories of the estimates  $\widehat{\x}(t)$ of the projected goal onto the collision-free local workspace, and the trajectories of the center of mass of the robot $\x_{\textsf{r}}(t)$. %The blue circle represents a particular configuration of the robot, where the collision-free local workspace $\mathcal L \mathcal F (\x_{\textsf{r}})$ is the purple polygon.}

In Figure~\ref{fig_dynamic_goal_trajectory}-right, we plot the metric defined in \eqref{eq_ate} for the algorithm with and without prediction; that is, structured and unstructured, respectively. Evidently, there is significant benefit using the structured algorithm.%: computing the SG metric on the mean tracking error over the simulation time, we observe a SG$\sim 20$.   

}

% !TEX root = revised_final.tex

%%%%%%%%%%%%%%%%%%%%%%%%%%%%%%%%%%%%%%%%%%%%%%%%%%%%%%%%%%%%%%%%%%%%%%%%%%%%%%%%
%%%%%%%%%%%%%%%%%%%%%%%%%%%%%%%%%%%%%%%%%%%%%%%%%%%%%%%%%%%%%%%%%%%%%%%%%%%%%%%%
%%%%%%%%%%%%%%%%% CONCLUSIONS
%%%%%%%%%%%%%%%%%%%%%%%%%%%%%%%%%%%%%%%%%%%%%%%%%%%%%%%%%%%%%%%%%%%%%%%%%%%%%%%%
%%%%%%%%%%%%%%%%%%%%%%%%%%%%%%%%%%%%%%%%%%%%%%%%%%%%%%%%%%%%%%%%%%%%%%%%%%%%%%%%

\section{Research outlook and future challenges}

Time-varying optimization is rapidly arising as an attractive algorithmic framework for today's fast-changing complex systems and world-size networks that entail heterogeneous and spatially distributed data streams. This article delineated the framework and underlined that structured algorithms can offer improved solutions to time-varying problems. \rev{In this section, a brief and certainly non-exhaustive list of the current challenges for structured and unstructured methods is outlined, with due implications in a  number of potential applications. }

\rev{
\smallskip

\noindent {\bf Wider classes of problems.} It has been already argued that unstructured methods generally require less functional assumptions than structured ones. For example, unstructured methods have been proposed for various non-strongly convex problems, as well nonconvex cost functions, where notions of dynamic regret can be used as performance indicators, see e.g.,~\cite{Shahrampour2018,Bernstein2019}, and e.g.,~\cite{Liu2018,Tang2018,Fattahi2019,Akhriev2020}. An attractive feature of time-varying nonconvex optimization algorithms is that they can be free of locally optimal trajectories. For structured methods, these classes of problems are largely unexplored, since e.g., underlying evolution models will have to be set-valued for non-strictly convex time-varying problems (because the solution trajectory is not unique). Interesting questions regarding bifurcations and merging of locally optimal trajectories, as well as the possibility of escaping isolated locally optimal trajectories naturally arise in this setting. A few efforts in this direction include  in~\cite{Guddat1990,Allgower1990,Dontchev2013}, but a comprehensive framework is lacking. A possible venue in this area could rely on piecewise linear continuation methods~\cite{Allgower1990}.  

\smallskip

\noindent {\bf Data-driven models.} Dynamic means of capturing the underlying optimization trajectory are now largely based on models, while in the current data streaming era, problems are often constructed in a data-driven fashion (e.g., via zero-order/bandit methods~\cite{Flaxman2005} or in a Bayesian setting~\cite{Srinivas2012}). Constructing and learning dynamic models for the optimization trajectory (for instance via historical data) is a largely unexplored territory, especially for structured methods-- where high-order smoothness is required for enhanced performance, in contrast with what typically (noisy) zero-order methods can provide. Unstructured methods can be found in~\cite{Slivkins2008, Shames2019, Bogunovic2016}. %{\color{red} Emiliano: here feel free to expand. [I think it is just fine]}

%\smallskip

%\noindent {\bf Different solution methods.} Most current structured methods rely on models based on Taylor expansions of optimality conditions. This eases the convergence analysis, but restrains their applicability.  

\smallskip

\noindent {\bf Distributed architectures.} Distributed methods to solve time-varying optimization problems (possibly involving large-scale networks) are key in many contemporary cyber-physical applications. Both structured and unstructured methods have been investigated~\cite{Rahili2015, Shahrampour2018, Gong2016, Xi2016a, Sun2017, Maros2019, Maros2017, Ling2013, Jakubiec2013, Zavlanos2013, Paper2, bastianello2020distributed}, but many challenges remain. As discussed in~\cite{SPM}, most distributed methods rely on diminishing step-size rules, which might not be an appropriate choice in time-varying settings when the algorithm runs indefinitely (as in e.g., video surveillance and monitoring of critical infrastructure). Another insight from~\cite{SPM,bastianello2020distributed} and~\cite{yuan2020can} is that the convergence behavior of distributed algorithms in the online setting is different relative to the batch case: traditional hierarchies in terms of convergence may be ``flipped,'' with the slowest algorithm in the static case being the fastest algorithm in the time-varying one. In addition, the notion of asynchronous updates assumes a more prominent position, inasmuch the network of computing nodes may have access to different evolution models, sample the optimization problem at different time steps, at different time-scales, or deliver solutions with different accuracy. All of this hinders standard analysis and it remains largely unexplored.  %{\color{red} All: feel free to expand, add more citations. [Emiliano: done]}

\smallskip

\noindent {\bf Feedback loop.} As we have seen in the analytical results presented here, under the assumptions provided, the time-varying algorithms converge to an error bound. Two key aspects are that \emph{(i)} the error bound can be arbitrarily big, if the algorithm converges arbitrarily slow, that is, if $\varrho$ is arbitrarily close to $1$; \emph{(ii)} the time-varying algorithms are considered separately, meaning the decision stream $\hat{\x}(t)$ does not influence the optimization problem at future times. Ensuring ``close-loop'' stability and performance, when the decision stream is fed back to the system is a mostly open challenge, and one can expect that arbitrarily slow algorithms cause lack of convergence. In this case, the very notion of ATE may be ill-defined or too hard to achieve, since typically the cost will be parametrized also on the approximated optimizer trajectory, and system-oriented notions of stability and robustness may be more appropriate. Some initial work can be found in~\cite{Hours2014,paternain2019prediction} in the context of MPC, yet this area remains largely open. %{\color{red} Santiago: feel free to expand on citations/ text. }
 
{\color{black} Another emerging research topic is the development of online structured and unstructured online algorithms that effectively act as feedback controlled dynamical systems. The main goal is to drive the output of a dynamical system  to solutions of time-varying optimization problems. Initial efforts toward unstructured online algorithms include~\cite{colombino2019online,zheng2019implicit}, where a Lyapunov analysis is also provided, and the more recent works in~\cite{nonhoff2019online,agarwal2019logarithmic} that provide a pertinent regret analysis).  
}
 
\smallskip

\noindent {\bf Interactive and Reinforcement Learning.} 
We close with potential links of the time-varying optimization tools outlined in this work with related contemporary thrusts on online convex optimization (OCO)~\cite{Shalev-Shwartz2012}, bandits~\cite{Agrawal2020}, and reinforcement learning (RL) that encompasses interactive decision making between agents and generally dynamic environments~\cite{Sutton2017}. At this stage, these links are active research thrusts that are pursued in diverse applications, such as allocation of network resources, secure mobile edge computing, and management of Internet-of-Things; see e.g., ~\cite{chen2019pieee,lcg2019tsp}, and references therein. Clearly, at the core of OCO, bandits, and RL are sequential solutions of optimization objectives that vary as the environment transitions across states and the agents take actions dynamically. These key elements prompt one to foresee that the time-varying tools overviewed in the present article can be fruitfully leveraged in interactive optimization. One key challenge to bear in mind in this direction is that the objective function in RL changes not only due to time-varying effects, but also due to actions fed back by the agent (learner). How to broaden the scope of algorithms presented here in such a wider context, constitutes an exciting open research direction.  
}

\bibliographystyle{ieeetr}
\bibliography{PaperCollection00}

\vspace*{-0.5cm}
\begin{IEEEbiography}[{\includegraphics[width=1in,height=1.25in,clip,keepaspectratio]{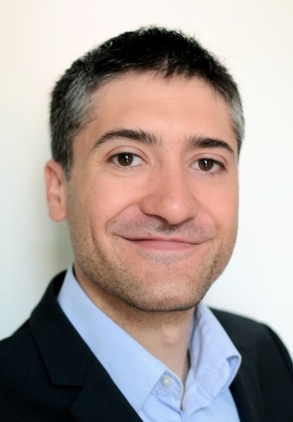}}]{Andrea~Simonetto} (M'12) is a research staff member in the optimization and control group of IBM Research Ireland, in Dublin.
He received his PhD in systems and control from Delft University of Technology, The Netherlands in 2012, and spent 3+1 years as postdoc, first in the signal processing group in the electrical engineering department in Delft, then in the applied mathematics department of the Universit\'e catholique de Louvain, in Belgium. %He was a visiting researcher at Carnegie Mellon University, University of Pennsylvania, and KTH, Sweden.
He joined IBM Research in February 2017.
His interests span optimization, control, and signal processing, with applications in smart energy, transportation, and personalized health.
\end{IEEEbiography}

\begin{IEEEbiography}[{\includegraphics[width=1in,height=1.25in,clip,keepaspectratio]{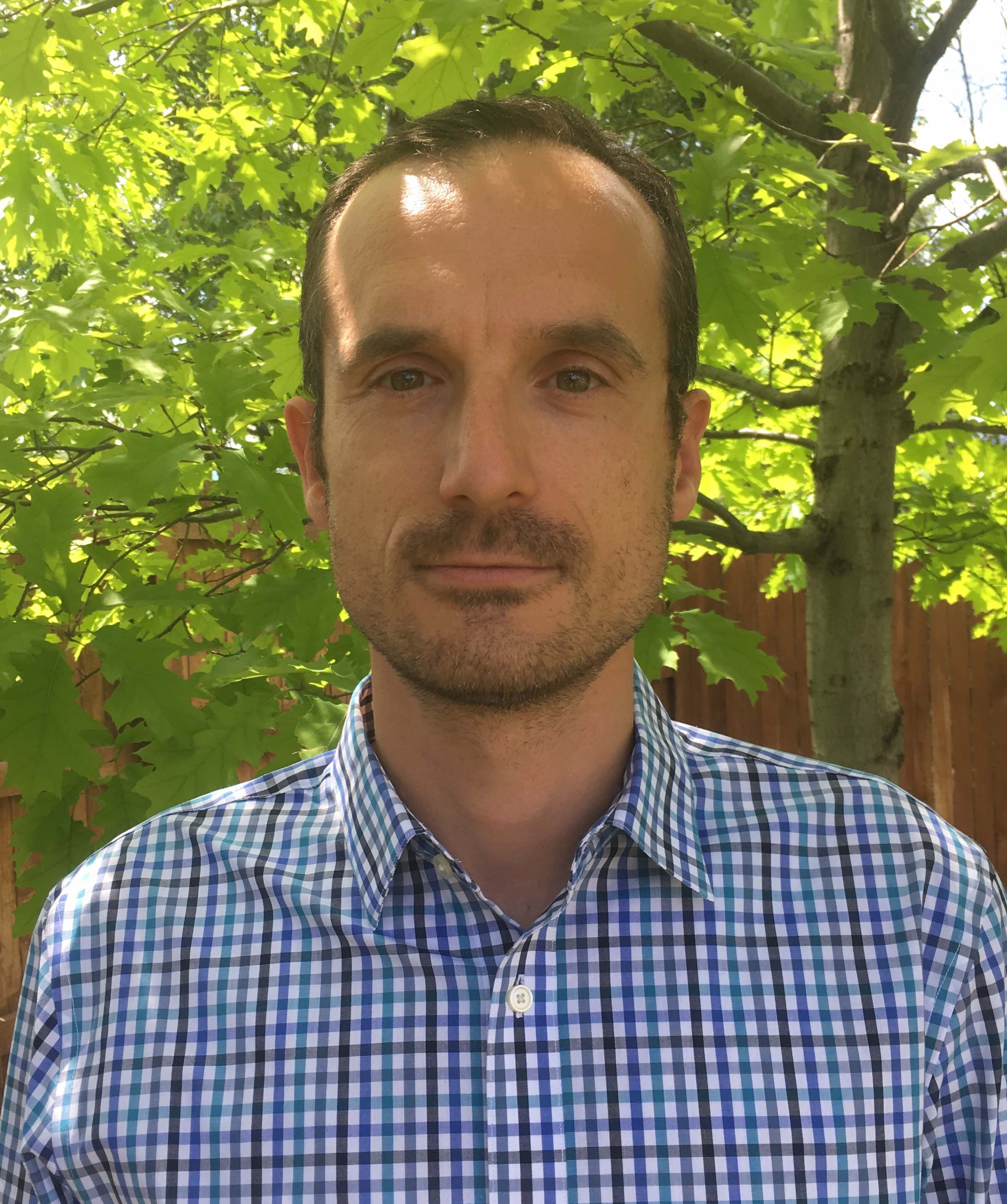}}]{Emiliano~Dall'Anese} (M'11) is an Assistant Professor within the Department of Electrical, Computer, and Energy Engineering at the University of Colorado Boulder. He received  the Ph.D. in Information Engineering from the Department of Information Engineering, University of Padova, Italy, in 2011. From January 2009 to September 2010, he was a visiting scholar at the Department of Electrical and Computer Engineering, University of Minnesota, USA. From January 2011 to November 2014, he was a Postdoctoral Associate at the Department of Electrical and Computer Engineering of the University of Minnesota, USA, and from December 2014 to July 2018 he was a Senior Researcher at the National Renewable Energy Laboratory, Golden, CO, USA. His research interests span the areas of  optimization,  control,  and signal processing,  with applications to networked systems and energy systems.  
\end{IEEEbiography}

\begin{IEEEbiography}[{\includegraphics[width=1in,height=1.25in,clip,keepaspectratio]{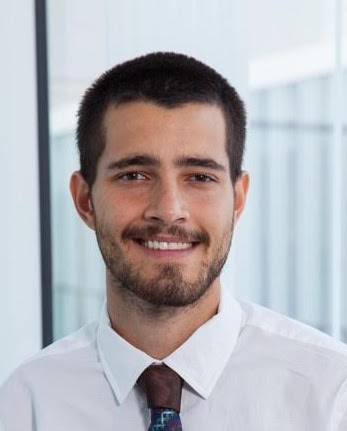}}]{Santiago~Paternain} (M'15) received the B.Sc. degree in electrical engineering from Universidad de la Rep\'ublica Oriental del Uruguay, Montevideo, Uruguay in 2012 and the M.Sc. in Statistics from the Wharton School in 2018 and the Ph.D. in electrical and systems engineering from the Department of Electrical and Systems Engineering, the University of Pennsylvania in 2018. He was the recipient of the 2017 CDC Best Student Paper Award.  His research interests include optimization and control of dynamical systems. 
\end{IEEEbiography}

\begin{IEEEbiography}[{\includegraphics[width=1in,height=1.25in,clip,keepaspectratio]{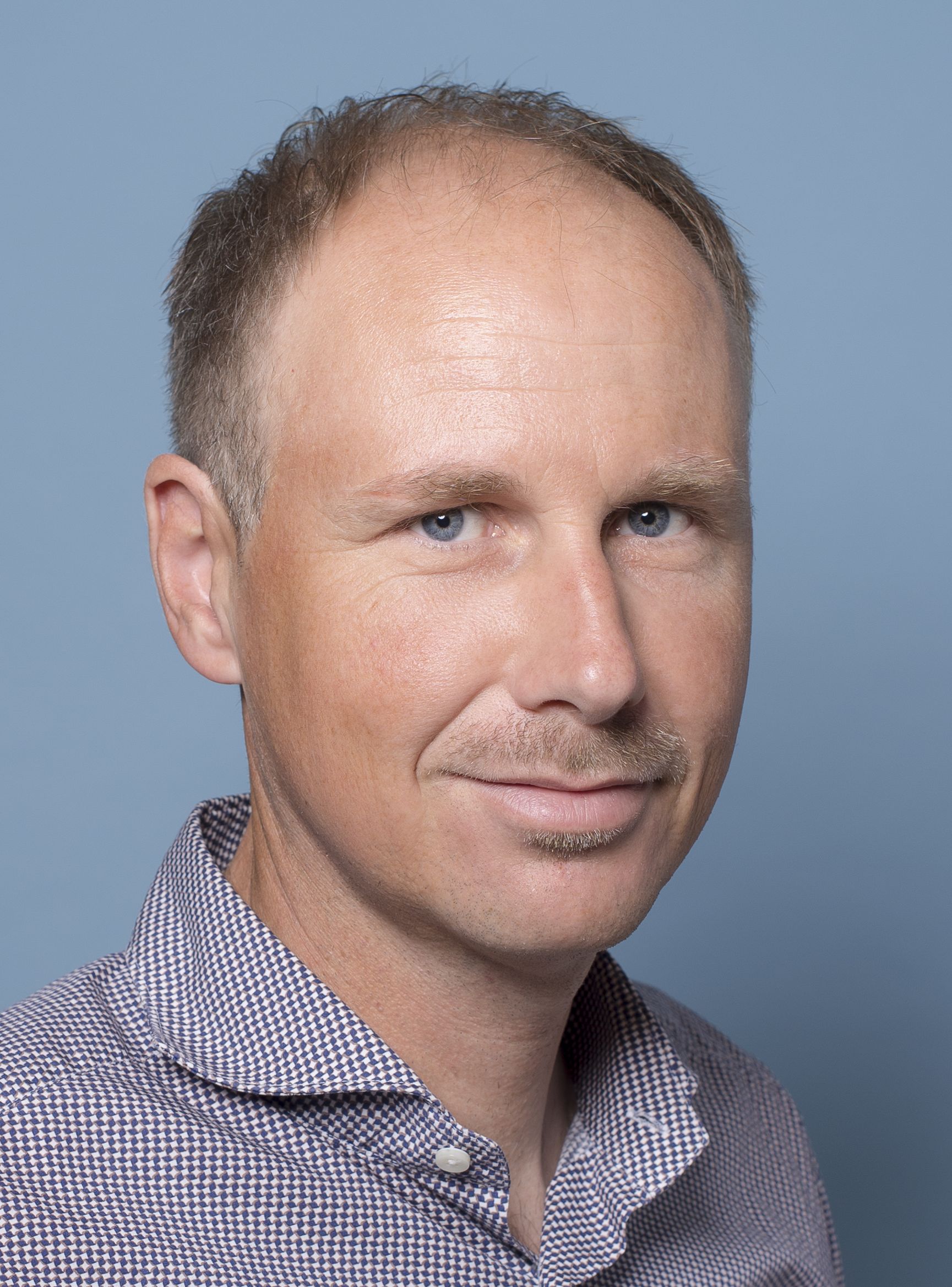}}]{Geert~Leus} (F'12) received the M.Sc. and Ph.D. degree in Electrical Engineering from the KU Leuven, Belgium, in June 1996 and May 2000, respectively. Geert Leus is now an ``Antoni van Leeuwenhoek'' Full Professor at the Faculty of Electrical Engineering, Mathematics and Computer Science of the Delft University of Technology, The Netherlands. His research interests are in the broad area of signal processing, with a specific focus on wireless communications, array processing, sensor networks, and graph signal processing. Geert Leus received a 2002 IEEE Signal Processing Society Young Author Best Paper Award and a 2005 IEEE Signal Processing Society Best Paper Award. He is a Fellow of the IEEE and a Fellow of EURASIP. Geert Leus was a Member-at-Large of the Board of Governors of the IEEE Signal Processing Society, the Chair of the IEEE Signal Processing for Communications and Networking Technical Committee, a Member of the IEEE Sensor Array and Multichannel Technical Committee, and the Editor in Chief of the EURASIP Journal on Advances in Signal Processing. He was also on the Editorial Boards of the IEEE Trans. on Signal Processing, the IEEE Trans. on Wireless Communications, the IEEE Signal Processing Letters, and the EURASIP Journal on Advances in Signal Processing. Currently, he is the Chair of the EURASIP Special Area Team on Signal Processing for Multisensor Systems, a Member of the IEEE Signal Processing Theory and Methods Tech. Committee, a Member of the IEEE Big Data Special Interest Group, an Associate Editor of Foundations and Trends in Signal Processing, and the Editor in Chief of EURASIP Signal Processing.  
\end{IEEEbiography}

\begin{IEEEbiography}[{\includegraphics[width=1in,height=1.25in,clip,keepaspectratio]{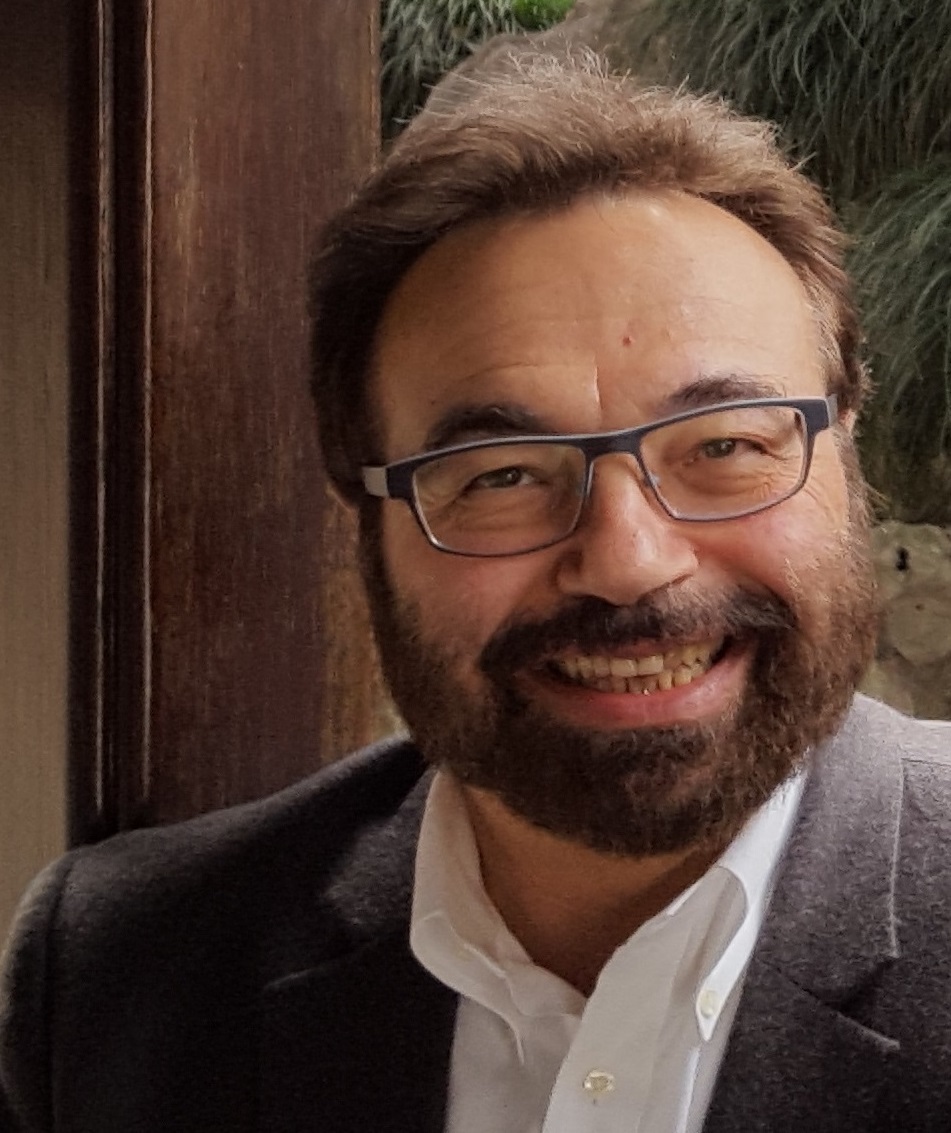}}]{Georgios~B.~Giannakis} (F'97) received his Diploma in Electrical
	Engr. from the Ntl. Tech. Univ. of Athens, Greece, 1981. From
	1982 to 1986 he was with the U. of Southern California (USC),
	where he received his MSc. in Electrical Engineering, 1983, MSc.
	in Mathematics, 1986, and Ph.D. in Electrical Engr., 1986. He
	was a faculty member with the U. of Virginia from 1987
	to 1998, and since 1999 he has been a professor with the U.
	of Minnesota, where he holds an ADC Endowed Chair, a University
	of Minnesota McKnight Presidential Chair in ECE, and serves as
	director of the Digital Technology Center. 
	His general interests span the areas of statistical learning,
	communications, and networking - subjects on which he has published
	more than 465 journal papers, 765 conference papers, 25 book
	chapters, two edited books and two research monographs. Current
	research focuses on Data Science, and Network Science with
	applications to the Internet of Things, and power networks with
	renewables. He is the (co-) inventor of 33 issued patents, and
	the (co-) recipient of 9 best journal paper awards from the IEEE
	Signal Processing (SP) and Communications Societies, including
	the G. Marconi Prize Paper Award in Wireless Communications. He also
	received the IEEE-SPS Nobert Wiener Society Award (2019); EURASIP's
	A. Papoulis Society Award (2020); Technical Achievement Awards
	from the IEEE-SPS (2000) and from EURASIP (2005); the IEEE ComSoc
	Education Award (2019); the G. W. Taylor Award for Distinguished
	Research from the University of Minnesota, and the IEEE Fourier
	Technical Field Award (2015). He is a Fellow of the National
	Academy of Inventors, the European Academy of Sciences, IEEE and
	EURASIP. He has served the IEEE in a number of posts, including
	that of a Distinguished Lecturer for the IEEE-SPS.

\end{IEEEbiography}

\end{document}